%% file: integral_differential_forms.tex
\documentclass{article}
\usepackage[utf8]{inputenc}
\usepackage{calligra,mathrsfs}
\usepackage{amssymb}
\usepackage{amsthm}
\usepackage{latexsym}
\usepackage{amsxtra}
\usepackage{fullpage}
\usepackage{amscd}
\usepackage{xypic}
\usepackage{verbatim}
\usepackage{enumerate}
\usepackage{url}
\usepackage[usenames,dvipsnames]{color}
\usepackage{tikz}
\usepackage{tikz-cd}
\usepackage{subfig}

\tikzset{
	solid node/.style={circle,draw,inner sep=1.2,fill=black},
	gray node/.style={circle,draw=gray,inner sep=1.2, fill = gray},
}

\usepackage{hyperref}

\newtheorem{thm}{Theorem}[section]
\newtheorem{lem}[thm]{Lemma}
\newtheorem{cor}[thm]{Corollary}
\newtheorem{prop}[thm]{Proposition}
\theoremstyle{definition}
\newtheorem{defn}[thm]{Definition}
\newtheorem{exa}[thm]{Example}
\newtheorem{rem}[thm]{Remark}
\newtheorem{ass}[thm]{Assumption}

\newtheorem{algorithm}[thm]{Algorithm}

\numberwithin{equation}{section}

\renewcommand{\em}{\sl}

\makeatletter
\renewcommand{\subsection}{\@startsection{subsection}{2}%
        {\z@}{-3.25ex plus -1ex minus-.2ex}{-1em}{\bf}}
\makeatother

\newcommand{\ZZ}{\mathbb{Z}}
\newcommand{\NN}{\mathbb{N}}
\newcommand{\QQ}{\mathbb{Q}}
\newcommand{\FF}{\mathbb{F}}
\newcommand{\RR}{\mathbb{R}}
\newcommand{\PP}{\mathbb{P}}

\renewcommand{\AA}{\mathbb{A}}

\newcommand{\OO}{\mathcal{O}}
\newcommand{\oo}{\mathfrak{o}}
\newcommand{\X}{\mathcal{X}}

\newcommand{\Spec}{{\rm Spec}}
\newcommand{\oomega}{\mbox{\Large $\omega$}}

\newcommand{\QQh}{\hat{\QQ}}
\newcommand{\Ab}{\bar{A}}

\newcommand{\p}{\mathfrak{p}}
\newcommand{\m}{\mathfrak{m}}

\newcommand{\vg}{{\mathbf{v}_0}}

\newcommand{\an}{{\rm an}}

\newcommand{\sing}{{\rm sing}}

\newcommand{\hor}{{\rm hor}}
\newcommand{\ord}{{\rm ord}}

\newcommand{\pred}{{\rm pred}}
\newcommand{\lcm}{{\rm lcm}}
\newcommand{\rad}{{\rm rad}}

\newcommand{\inj}{\hookrightarrow}
\newcommand{\iso}{\stackrel{\sim}{\to}}

\newcommand{\HOM}{\mathscr{H}\text{\kern -3pt {\calligra\large om}}\,}

\newcommand{\gen}[1]{\left\langle #1 \right\rangle}
\newcommand{\abs}[1]{\left\lvert #1 \right\rvert}
\newcommand{\norm}[1]{\left\lVert #1 \right\rVert}

\begin{document}

\title{Integral differential forms for superelliptic curves}

\author{Sabrina Kunzweiler and Stefan Wewers}

\date{}
\maketitle

\begin{abstract}
	Given a superelliptic curve $Y_K:\; y^n=f(x)$	
	over a local field $K$, we describe the theoretical background and an implementation of a new algorithm for computing the $\oo_K$-lattice of integral differential forms on $Y_K$. We build on the results of \cite{ObusWewers} which describe arbitrary regular models of the projective line using only valuations. One novelty of our approach is that we construct an $\oo_K$-model of $Y_K$ with only rational singularities, but which may not be regular.    
	
	\vspace{2ex}\noindent 
	2020 {\em Mathematics Subject Classification}. Primary 11G20. 
	Secondary: 14G10.
\end{abstract}

\section{Introduction} 

\subsection{Regular models, rational singularities, and integral differential forms}

Let $K$ be a local field with ring of integers $\oo_K$ and $Y_K$ a smooth projective curve of genus $g\geq 1$ over $K$. The group of global sections of the sheaf $\Omega_{Y_K/K}$ of differential $1$-forms,
\[
   M_K := H^0(Y_K,\Omega_{Y_K/K}),
\] 
is a $K$-vector space of dimension $g$. It is equipped with a natural $\oo_K$-lattice $M\subset M_K$, the lattice of {\em integral differential forms}. We recall the definition of $M$.

By a {\em model} of $Y_K$ we mean a proper, flat and normal relative curve $Y\to S:=\Spec\,\oo_K$ with generic fiber $Y_K$. On such a model there exists a natural extension $\oomega_{Y/S}$ of $\Omega_{Y_K/K}$ called the {\em canonical sheaf}. It is a coherent $\OO_Y$-sheaf such that $\oomega_{Y/S}|_{Y_K} = \Omega_{Y_K/K}$. See e.g.\ \cite[\S 6.4]{liu2002algebraic} and \cite{morrow2010explicit}. The group of global sections,
\[
   M_Y := H^0(Y,\oomega_{Y/S}),
\]
is a full $\oo_K$-lattice of $M_K$. 

Let $Y$, $Y'$ be models of $Y_K$. We say that $Y'$ {\em dominates} $Y$ if the identity on $Y_K$ extends (uniquely) to a birational
morphism $f:Y'\to Y$. If this is the case then there is a natural embedding of sheaves
\[
     f_*\oomega_{Y'/S} \inj \oomega_{Y/S},
\]
resulting in an inclusion
\[
   M_{Y'} \subset M_Y.
\]

A point $y\in Y$ is a {\em regular point} if $\OO_{Y,y}$ is a regular local ring. Otherwise, $y$ is called a {\em singularity}. If every point of $Y$ is regular, we say that $Y$ is a {\em regular model}. If $Y$ is an arbitrary model of $Y_K$, a {\em desingularization} of $Y$ is a regular model $Y'$ dominating $Y$. It is known that every model $Y$ has a desingularization (\cite{Lipman78}).

We say that a point $y\in Y$ is a {\em rational singularity} if for every desingularization $f:Y'\to Y$ we have
\[
   (f_*\oomega_{Y'/S})_y = (\oomega_{Y/S})_y.
\]
So if $Y$ has only rational singularities then 
\[
   f_*\oomega_{Y'/S} = \oomega_{Y/S}
\]
holds for every desingularization $f:Y'\to Y$.

\begin{prop}
  A regular point $y\in Y$ is a rational singularity.
\end{prop}

\begin{proof}
The canonical sheaf $\oomega_{Y/S}$ is the {\em dualizing sheaf} for the proper morphism $Y\to S$, see \cite[\S 6.4.3]{liu2002algebraic} and \cite{morrow2010explicit}. By definition, $y$ is a rational singularity if and only if for every desingularization $f:Y'\to Y$ we have
\[
    (f_*\oomega_{Y'/S})_y = (\oomega_{Y/S})_y.
\]
By duality,
 this holds if and only if
\[
   (R^1f_*\OO_{Y'})_y = 0.
\]
(The latter condition is actually the usual definition of the term {\em rational singularity}). The statement of the proposition follows now from \cite[Proposition 1.2]{lipman1969rational}.
\end{proof}

\begin{cor}
	\label{cor:lattice_rationalsingularities}
  	A model $Y$ of $Y_K$ has only rational singularities if and only if the $\oo_K$-lattice
  	\[
  	    M=M_Y:=H^0(Y,\oomega_{Y/S})\subset M_K
  	\]
  	is minimal, among all lattices of this form. Moreover, this is the case if $Y$ is regular.
\end{cor}

\begin{defn}
 \label{def:integral_differential_forms}
  The lattice of {\em integral differential forms} for the curve $Y_K$ is the lattice $M_Y\subset M_K$, where $Y$ is any model of $Y_K$ with at most rational singularities.	
\end{defn}

It is an important problem in computational arithmetic geometry to compute the lattice $M$ explicitly for a given curve $Y_K$. One reason for this is that $M$ is equal to the group of global invariant differential forms on the N\'eron model of the Jacobian of $Y_K$ (see e.g.\ \cite[\S 3.2]{vanBommel2019numerical}). Computing (the top exterior power of) the latter is an important step towards the numerical verification of the Birch and Swinnerton-Dyer conjecture for the curve $Y_K$ (see e.g.\ \cite{flynn2001empirical}, \cite{vanBommel2019numerical}).

\vspace{2ex}
For hyperelliptic curves with semistable reduction, there exists a very compact formula that allows to determine the top exterior power of $M$, see \cite{kausz1999discriminant} and \cite{kunzweiler2020differential}. In a more general setting, the strategy of choice for computing the lattice $M$ seems to have been to compute a regular model of $Y_K$, using one of the existing implementations, see e.g.~\cite{flynn2001empirical}, \cite{vanBommel2019numerical}, \cite{dokchitser2021models}, \cite{muselli2020models}.
In \cite{flynn2001empirical} and \cite{vanBommel2019numerical}, the implementation of regular models by Steve Donnelly in Magma (\cite{magma}) is used. The algorithm behind this implementation (successively blowing up a model in the singular points) is completely general. However, the current implementation requires some strong restrictions on the curve $Y_K$ (essentially it needs to be a smooth projective plane curve). Tim Dokchitser's Magma script \cite{dokchitser_regular_models}, which is based on the preprint \cite{dokchitser2021models} requires that the curve $Y_K$ is given as a smooth affine curve in $\mathbb{G}_m^2$ and satisfies an extra `genericity' condition. It is not clear to us how restrictive this assumption is. The method presented in \cite{muselli2020models} combines ideas both from \cite{dokchitser2021models} and \cite{MMDD}. So far it has only been worked out in the case of hyperelliptic curves. 

In the present article, we suggest an alternative approach for computing the lattice of integral differential forms. It applies to some new cases not covered by the implementations mentioned above. It is based on the following two observations:
\begin{enumerate}[(a)]
\item 
  It is often easier to find a model $Y$ with at most rational singularities (which may not be regular).
\item 
  In order to compute the lattice $M=M_Y$, it  suffices to find explicit equations for the model $Y$ in some \'etale neighborhood of the generic points of the special fiber. This can also be a lot easier than to find equations for a Zariski covering of the whole model. 
\end{enumerate}	 
Both points are rather obvious, but we hope to show their usefulness by working out a special case in detail and providing a concrete implementation.

\subsection{Superelliptic curves}

We continue with the notation used before. The curve $Y_K$ is called {\em superelliptic} if it is the smooth projective model of a plane affine curve given by an equation of the form
\begin{equation} \label{eq:superelliptic1}
    y^n = f(x),
\end{equation}
where $n\geq 2$ and $f\in K[x]$ is a polynomial of the form
\[
     f = \prod_{i=1}^r f_i^{m_i},
\]
with $f_i\in K[x]$ irreducible, $1\leq m_i<n$ and $\sum_i\deg f_i\geq 3$. For instance, if $n=2$ then $Y_K$ is a {\em hyperelliptic curve}. The crucial assumption we impose on $Y_K$ is:

\begin{ass} \label{ass:tame}
  The exponent $n$ in \eqref{eq:superelliptic1} is invertible in the ring $\oo_K$.
\end{ass}

Under this assumption, we can construct a model $Y$ of $Y_K$ with at most rational singularities as follows. The equation \eqref{eq:superelliptic1} presents $Y_K$ as a cyclic cover of the projective line, of degree $n$,
\[
    \phi_K:Y_K \to X_K := \PP^1_K.
\]
Let $D_K\subset X_K$ denote the {\em branch locus} of this cover (which is equal to the divisor of zeroes of $\rad(f)=\prod_i f_i$, possibly joined by the point $\infty$). Let $X$ be a model of $X_K$. We denote by $D$ the Weil divisor
\begin{equation} \label{eq:divisor}
     D = D^{\rm vert} \cup D^{\rm hor}
\end{equation}
on $X$, where $D^{\rm vert}:=(X\otimes k)_{\rm red}$ is the reduced special fiber of $X$ and $D^{\rm hor}$ is the Zariski closure of $D_K\subset X_K$ inside $X$. Moreover, let $Y$ denote the normalization of $X$ inside the function field of $Y_K$. This is the unique model of $Y_K$ such that the cover $\phi_K$ extends to a finite map 
\[
    \phi:Y\to X.
\]

\begin{thm} \label{thm:tame}
  Assume the following:
  \begin{enumerate}[(a)]
  \item 
    $n$ is invertible in $\oo_K$ (Assumption \ref{ass:tame}).
  \item 
    The model $X$ is regular.
  \item 
    The divisor $D$ is a normal crossing divisor.    	
  \end{enumerate}	
  Let $U:=X\backslash D^\sing$ denote the complement of the singular points of the divisor $D$, and $V:=\phi^{-1}(U)$. Then $V$ is regular, and the points in $Y\backslash V$ are rational singularities.
  
  In particular, the model $Y$ has at most rational singularities.
\end{thm}

\begin{proof}
It follows from the Assumptions (a)-(c) and a version of  {\em Abhyankar's Theorem} (\cite[Corollary 2.3.4]{GrothendieckMurre}) that the map $\phi:Y\to X$ is
\begin{enumerate}[(i)]
\item
  \'etale over $X\backslash D$,
\item 
  a {\em Kummer cover} in an \'etale neighborhood of a smooth point of $D$, and
\item 
  a   {\em generalized Kummer cover}, as defined in \cite[\S 1]{GrothendieckMurre}, in an \'etale neighborhood of a singular point of $D$.
\end{enumerate}
Here (ii) means that the cover is given \'etale locally by an equation $T^m=t$, where $t$ is a local equation for the divisor $D$ and $m\mid n$. It follows immediately that $V$ is regular. 

Let $x\in D^\sing$ be a singular point of $D$. Then in an \'etale neighborhood of $x$, the divisor $D$ is given by an equation $t_1t_2=0$, where $(t_1,t_2)$ is a system of parameters for the regular point $x\in X$. Then (iii) means that there exist (after replacing $X$ by a sufficiently small \'etale neighborhood of $x$ and $Y$ by a connected component of the inverse image)
\begin{itemize}
\item
	positive integers $n_1,n_2$, invertible on $X$,
\item 
	roots of unity $\zeta_{n_1},\zeta_{n_2}\in\OO_X$, of order $n_1, n_2$, 
\item 
	a cyclic subgroup $G=\gen{g}\subset \ZZ/n_1\ZZ\times\ZZ/n_2\ZZ$, and 
\item 
	a finite covering of $X$-schemes
	\[
	  Z:= X[T_1,T_2 \mid T_1^{n_1}=t_1,\,T_2^{n_2}=t_2] \to Y
	\]  
\end{itemize}
such that the following holds. The generator $g=(k_1,k_2)\in \ZZ/n_1\ZZ\times\ZZ/n_2\ZZ$ of $G$ acts on $Z/X$ via
\[
g^*T_1=\zeta_{n_1}^{k_1}\cdot T_1, \quad 
g^*T_2=\zeta_{n_2}^{k_2}\cdot T_2.
\]
Moreover, $Y=Z/G$.

It is clear that $Z$ is regular. Therefore, the unique point $y\in Y$ lying over $x$ is a {\em tame quotient singularity}. By the proof of the main result of \cite{boutot1987singularites}, tame quotient singularities are rational singularities\footnote{In \cite{boutot1987singularites} all schemes are assumed to be of finite type over a field of characteristic zero, or a localization of such schemes. But the proof goes through essentially unchanged in our situation. 
}. This concludes the proof of the theorem.
\end{proof}

\begin{rem}
  The statement of Theorem \ref{thm:tame} applies more generally to covers of curves $\phi_K:Y_X\to X_K$ with monodromy group of order prime to the residue characteristic of $K$. For instance, if $Y_K$ is a smooth plane quartic that admits an action of the Klein $4$-group $V$, then $Y_K$ may be realized as a cover of the projective line with monodromy group $V$ and we have to assume that the residue characteristic of $K$ is odd. A classification of such curves and their stable models is given in \cite{bouw2020reduction}. The methods developed in our work also apply to this class of curves. An explicit example is worked out in \cite[Example 4.34]{kunzweiler2021thesis}.
  
  More generally if $Y_K$ is any smooth plane quartic and $\phi_K$ the projection to one of the coordinate axes, then Theorem \ref{thm:tame} applies if the residue characteristic of $K$ is $\geq 5$. 
\end{rem}   

\begin{rem}
  The model $Y$ from Theorem \ref{thm:tame} will in general not be regular, but the singularities that appear are rather nice. Apart from being rational, they are {\em toric singularities} in the sense of \cite{kato1994toric}, and are easy to resolve, see \cite[ \S 4.4.2]{halle2016neron}. We plan to extend our implementation to also compute regular models. 
\end{rem}	

\subsection{Models and valuations}

In order to use Theorem \ref{thm:tame}, we need to be able to find a regular model $X$ of $X_K=\PP^1_K$ such that the divisor $D$ (which depends on $D_K$ and the model $X$) is a normal crossing divisor. Typically, taking the obvious model $X:=\PP^1_{\oo_K}$ (which is smooth over $\oo_K$ and hence regular) won't work, because the horizontal part $D^{\rm hor}$ of $D$ (the closure of $D_K$) may not be regular, and may not intersect the vertical component $D^{\rm vert}=\PP^1_k$ transversally. 

Our approach for constructing a suitable model $X$ is an extension of the methods introduced in \cite{RuethThesis} and \cite{ObusWewers}. The first observation is that a model $X$ of $X_K$ is uniquely determined by the finite set $V(X)$ of discrete valuations on the function field $F_X=K(x)$ corresponding to irreducible components of the special fiber $X_s$ of $X$ (see \cite[\S 3]{RuethThesis}). For instance, the model $X=\PP^1_{\oo_K}$ of $X_K=\PP^1_K$ corresponds to the set $V(X)=\{\vg\}$, where $\vg$ is the {\em Gauss valuation} on $K(x)$ with respect to the parameter $x$, extending the valuation $v_K$.
For an arbitrary model $X$, a valuation $v\in V(X)$  corresponding to an irreducible component $E_v\subset X$ of the special fiber can be described very explicitly, using the results of MacLane  (\cite{MacLane36}). This description is equivalent to giving an explicit equation for the scheme $X$ which is valid in a Zariski neighborhood of the generic point of $E_v$. 

It is a very interesting problem to find explicit equations for $X$ in the neighborhood of an arbitrary point, knowing only the set $V(X)$. While this problem has not been solved in general, one can often extract all the information from the set $V(X)$ that one needs. For instance, in \cite[\S 7]{ObusWewers} it is shown how to find, given an arbitrary model $X$ of $X_K=\PP^1_K$, an enlargement $V'\supset V(X)$ such that the model $X'$ with $V(X')=V'$ is a desingularization $X'\to X$ of $X$. In \S \ref{sec:regular_models_of_projective_line} of the present article, this result is extended to a solution of the problem posed by Theorem \ref{thm:tame} and formulated at the beginning of this subsection. The essential new ingredient is to allow the set $V(X)$ to contain certain {\em infinite pseudovaluations} which correspond to the irreducible horizontal components of the divisor $D$ (i.e.\ to the points of $D_K$). We remark that a very similar result is obtained in \cite{ObusSrinivasan}, see also \cite{obus2022explicit}.

Once we have found a model $X$ satisfying the Conditions (b) and (c) of Theorem \ref{thm:tame}, we easily obtain the set $V(Y)$ of valuations corresponding to the vertical components of the induced model $Y$ of $Y_K$. In fact, $V(Y)$ is just the set of extensions of the valuations in $V(X)$ to the function field of $Y_K$. Since the extension $F_Y/F_X$ is given by the simple equation \eqref{eq:superelliptic1}, it is very easy to compute $V(Y)$\footnote{The methods from \cite{MacLane36} and \cite{RuethThesis} apply to arbitrary extensions of function fields.}

Finally, we come back to our original goal of computing the lattice $M$ of integral differential forms. Using the fact that the canonical sheaf $\oomega_{Y/S}$ is {\em divisorial},
 we have
\begin{equation} \label{eq:superelliptic2}
   M = \{ \omega\in M_K \mid v(\omega)\geq 0 \;\forall v\in V(Y)\}.
\end{equation}
Here $v(\omega)\in\ZZ$ denotes the order of vanishing of $\omega$ along the vertical component $E_v\subset Y$ corresponding to $v$, as a rational section of the canonical sheaf $\oomega_{Y/S}$. In \S \ref{sec:global_sections} of the present article we show how to use \eqref{eq:superelliptic2} to compute an explicit integral basis of $M$. 

\vspace{2ex}
We have implemented the method described above for computing a basis of integral differential forms for a superelliptic curve \eqref{eq:superelliptic1} satisfying Assumption \ref{ass:tame}. Our implementation is written in Python/Sage (\cite{sagemath}) and builds upon the Sage toolbox \href{https://github.com/MCLF/mclf}{MCLF: Models of Curves over Local Fields} (\cite{mclf}).

\vspace{2ex}
The article is organized as follows. In \S \ref{sec:inductive_valuations} we recall the theory of {\em inductive valuations} from \cite{MacLane36} and prove some technical results used later. In \S \ref{sec:regular_models_of_projective_line} we extend the results of \cite{ObusWewers} and develop an algorithm for computing a regular model of the projective line $X$ such that the divisor $D$ \eqref{eq:divisor} is a normal crossing divisor.
In \S \ref{sec:global_sections} we describe our main algorithm for computing the lattice of integral differential forms. Finally, in \S \ref{sec:examples} we illustrate our results by working out two examples in some detail. 

The results of this manuscript have been published in the first author's PhD thesis at the University of Ulm \cite{kunzweiler2021thesis}.

\subsection*{Acknowledgments} We thank Raymond van Bommel, Andrew Obus and the anonymous referee for valuable comments that helped to improve our manuscript.

\input{inductive_valuations}

\input{regular_models_of_P1}

\input{global_sections}

\input{examples_new}

\bibliographystyle{amsplain}
\bibliography{references}

\vspace{5ex}\noindent {\small Sabrina Kunzweiler\\  Inria Centre de l'universit\'e de Bordeaux\\ 351, cours de la Libération\\ 33405 Talence\\ {\tt
		sabrina.kunzweiler@math.u-bordeaux.fr}
	
\vspace{5ex}\noindent {\small Stefan
	Wewers\\ Institut f\"ur Algebra und Zahlentheorie\\ Universit\"at
	Ulm\\ Helmholtzstr.\ 18\\ 89081 Ulm\\ {
	\tt stefan.wewers@uni-ulm.de}}

\end{document}

%% file: inductive_valuations.tex

\section{Inductive valuations} \label{sec:inductive_valuations}

In this section we recall certain results on (pseudo)valuations on a polynomial ring which are originally due to MacLane (\cite{MacLane36}). These results have been known and used a lot by experts in valuation theory. In recent years, they also have found many applications in algorithmic number theory (see e.g.\ \cite{Fernandez_residual}) and arithmetic geometry (see e.g.\ \cite{RuethThesis}, \cite{ObusWewers}, \cite{ObusSrinivasan}).

In this and the next section, there is quite some overlap with the articles \cite{ObusWewers} and \cite{ObusSrinivasan}, the reason being that our main assumptions and also our focus are somewhat different. For instance, in this article we can't afford to assume that the residue field $k$ of our ground field $K$ is algebraically closed, a very helpful assumption made in \cite{ObusWewers} and \cite{ObusSrinivasan}. Another reason for being a bit repetitive is that we want to treat classical valuations and certain {\em pseudovaluations} on an equal footing.

\subsection{MacLane pseudovaluations} \label{subsec:maclane}

Throughout, $K$ denotes a field which is complete with respect to a normalized discrete valuation
\[
   v_K: K\to\QQh:=\QQ\cup\{\infty\}
\]
({\em normalized} means that $v_K(K^\times) = \ZZ$). We  denote by $\oo_K$ the valuation ring, $\pi_K$ an arbitrary uniformizer and $k:=\oo_K/(\pi_K)$ the residue field of $v_K$.

Let $K[x]$ be the polynomial ring over $K$ in some unknown $x$. A polynomial $f=\sum_i a_ix^i$ is said to be {\em integral} if $v_K(a_i)\geq 0$ for all $i$.

\begin{defn} \label{def:maclane_valuation}
  A map
  \[
      v:K[x]\to \QQh
  \]
  is called a {\em pseudovaluation} if $v(1) = 0$ and 
  \begin{enumerate}[(a)]
  \item
     $v(fg) = v(f) + v(g)$,
  \item
    $v(f+g) \geq \min(v(f),v(g))$, for all $f,g\in K[x]$. \end{enumerate}
 Here we use the usual conventions regarding the symbol $\infty$.
 Conditions (a) and (b) show  that  $I_v:=v^{-1}(\infty)\lhd K[x]$ is a prime ideal and that $v$ induces a  valuation
 \[
    \bar{v}:K[x]/I_v \to \QQh.
 \]
 If $I_v=(0)$ then $v=\bar{v}$ is simply a valuation in the usual sense. It extends uniquely to a valuation on the fraction field $K(x)$ of $K[x]$, which we will also denote by $v$. We write $k(v)$ for the residue field of $v$.

 If $I_v\neq (0)$ then we call $v$ an {\em infinite pseudovaluation}. If this is the case, $I_v\lhd K[x]$ is a maximal ideal and $K_v:=K[x]/I_v$ is a finite field extension of $K$. We also write $k(v)$ for the residue field of the induced valuation $\bar{v}$ on $K_v$.

 Assume that (a) and (b) hold. Then $v$ is called a {\em MacLane pseudovaluation} if moreover the following holds:
 \begin{enumerate}[(a)]\setcounter{enumi}{2}
  \item
    $v|_K = v_K$,
  \item
    $v(x)\geq 0$,
  \item
    $v$ is discrete, i.e.\ $v(K[x])\cap\QQ = \frac{1}{e_v}\ZZ$, for some positive integer $e_v$,
  \item
    if $I_v=(0)$, then the extension $k(v)/k$ has transcendence degree one (note that $k\subset k(v)$ because of Condition (c)).
  \end{enumerate}

  We write
  \[
      V(K[x])^* = V(K[x]) \;\dot{\cup}\; V(K[x])_\infty
  \]
  for the set of all MacLane pseudovaluations, given by the disjoint union of the subset of all MacLane valuations, $V(K[x])$, and the subset of the infinite MacLane pseudovaluations, $V(K[x])_\infty$.
  We define a partial order $\leq$ on $V(K[x])^*$ by
  \[
  v\leq v' \quad:\Leftrightarrow\quad
  v(f) \leq v'(f)\;\;\forall\, f\in K[x].
  \]
\end{defn}

\begin{exa} \label{exa:gauss-valuation}
  The {\em Gauss valuation} $\vg:K[x]\to\QQh$ is defined by
  \[
     \vg(\sum_{i=0}^d a_ix^i) := \min_i v_K(a_i),
  \]
  for polynomials $f=\sum_i a_ix^i\in K[x]$. One checks that $v$ is a MacLane valuation with residue field $k(\vg)=k(\bar{x})$ (where $\bar{x}$ is the image of $x$ in $k(\vg)$).

  It is immediately clear from the above definitions that the Gauss valuation $\vg$ is the least element of $V(K[x])^*$. Also, $f\in K[x]$ is integral if and only if $\vg(f)\geq 0$ if and only if $v(f)\geq 0$ for all $v\in V(K[x])$.
\end{exa}

\begin{rem} \label{rem:infinite_pseudovaluations}
  Let $v\in V(K[x])_\infty$ be an infinite pseudovaluation. Then $I_v=(g)$, for a unique monic, integral and irreducible polynomial $g$. Since $K_v:=K[x]/(g)$ is a finite field extension of $K$ and $K$ is complete with respect to $v_K$, the induced valuation $\bar{v}$ is the unique extension of $v_K$ to $K_v$. It follows that $v$ is uniquely determined by $g$. Conversely, if $g\in K[x]$ is monic, irreducible and integral, then there exists a unique infinite MacLane pseudovaluation $v$ such that $v(g)=\infty$.

  It also follows that an infinite pseudovaluation is a maximal element of $V(K[x])$ (we will see later that they are exactly the maximal elements).
\end{rem}

\begin{rem} \label{rem:berkovich}
  To $v\in V(K[x])^*$ we can associate the {\em multiplicative seminorm}
  \[
      \norm{\cdot}_v: K[x]\to\RR_{\geq 0}, \qquad \norm{f}_v:=q^{-v(f)},
  \]
  where $q>1$ is fixed. By Condition (c), $\norm{\cdot}_v$ is an extension of the nonarchimedian absolute value $\norm{\cdot}$ on $K$. Therefore, we may regard $V(K[x])^*$ as a subset of the analytic $K$-space $(\AA_K^1)^{\rm an}$, as defined in \cite{Berkovich}. Condition (d) means that $\norm{\cdot}_v$ lies on the {\em closed unit disk}, and Condition (e) and (f) mean that $\norm{\cdot}_v$ is either a {\em point of type I} (if $v$ is an infinite pseudovaluation), or a {\em point of type II} (if $v$ is a MacLane valuation). See \cite[\S 4.2]{Berkovich}.

  Some of the results in this article could probably be expressed more elegantly using the language of nonarchimedian analytic spaces. However, since the core of this work deals very explicitly with MacLane valuations as defined above, we decided to phrase everything in this language.
\end{rem}

\subsection{Key polynomials and augmentation} \label{subsec:augmentations}
A MacLane valuation $v \in V(K[x])$ can be augmented with a key polynomial $\phi \in K[x]$ to obtain a valuation $v' \geq v$. This process is explained in the following.

\begin{defn} \label{def:key_polynomial}
Let $v\in V(K[x])$ be a MacLane valuation and $f,g,h\in K[x]$.
\begin{enumerate}
\item
  We say that $f,g$ are {\em $v$-equivalent} (written as $f\sim_v g$) if $v(f-g)>v(f)=v(g)$.
\item
  The polynomial $f$ is said to {\em $v$-divide} $g$ (written as $f|_v\, g$) if $g\sim_v fh$, for some $h\in K[x]$. It is called {\em $v$-irreducible} if $f|_v\, gh$ implies $f|_v\, g$ or $f|_v\, h$, and {\em $v$-minimal} if $f|_v\,g$ implies $\deg(f)\leq\deg(g)$.
\item
  A polynomial $\phi\in K[x]$ is called a {\em key polynomial} for $v$ if $\phi$ is monic, integral, {\em $v$-irreducible} and {\em $v$-minimal}.
\end{enumerate}
See \cite[\S 1]{MacLane36}, or \cite[\S 4.1.1]{RuethThesis}.
\end{defn}

\begin{exa}
	Let  $\vg:\QQ_3[x]\to\QQh$ be the Gauss valuation (Example \ref{exa:gauss-valuation}). It is easy to see that all monic linear polynomials are key polynomials for $\vg$.

	The  polynomial $f = x^2+x+1 \in \QQ_3[x]$ is {\em not} a key polynomial for $\vg$ since it is not $\vg$-irreducible. To see this, note that
	\[
	0 = \vg((x+2)^2) = \vg(f) < \vg((x+2)^2 - f) = \vg(3x+3) = 1.
	\]

	On the other hand, the quadratic polynomial $g = x^2+1$ is a key polynomial for $\vg$. In general, one can show that any monic polynomial $\phi$ of positive degree is a key polynomial for $\vg$ if and only if $\vg(\phi) = 0$ and the reduction $\bar{\phi} \in \FF_3[x]$ is irreducible, see \cite[Lemma 4.8]{RuethThesis}.

\end{exa}

\begin{rem} \label{rem:key_is_irreducible}
One immediate consequence of these definitions is that any key polynomial $\phi$ for $v$ is irreducible (as an element of $K[x]$).
\end{rem}

Let $\phi$ be a key polynomial for $v\in V(K[x])$ and $\lambda\in\QQh$ such that $\lambda\geq v(\phi)$. Then one can define a new MacLane pseudovaluation
\[
v'=[v, v'(\phi)=\lambda] \in V(K[x])^*
\]
as follows. Any polynomial $f\in K[x]$ can be uniquely written as
\[
   f = \sum_i f_i \phi^i,
\]
with $\deg f_i<\deg \phi$. Then
\[
   v'(f) := \min_i (v(f_i) + i\cdot\lambda).
\]
See \cite[Definition 4.9]{RuethThesis}.
We call $v'$ the {\em augmentation} of $v$ with respect to $(\phi,\lambda)$. It holds that $v\leq v'$ and $v<v'$ if and only if $v(\phi)<\lambda$, see \cite[\S 5]{MacLane36} . In the latter case we say that $v'$ is a {\em proper augmentation}\footnote{In \cite{MacLane36} and \cite{RuethThesis}, only proper augmentations are called augmentations. Here we follow the convention introduced in \cite{ObusWewers}.} of $v$. Also, $v'$ is a MacLane valuation if and only if $\lambda\neq \infty$.

\begin{rem} \label{rem:augmented_valuations}
  If $v'=[v,v'(\phi)=\lambda]$ is an augmentation of $v$, then $\lambda$ and $\deg(\phi)$ only depend on $v$ and $v'$, but $\phi$ does not. More precisely, if $w=[v,w(\psi)=\mu]$ is another augmentation of $v$ then $v'=w$ if and only if $\lambda=\mu$, $\deg(\phi)=\deg(\psi)$ and
  \[
       v(\phi-\psi) \geq \lambda.
  \]
  See \cite[Theorem 15.3]{MacLane36} or \cite[Theorem 4.33]{RuethThesis}.
\end{rem}

\subsection{Representing MacLane pseudovaluations as inductive valuations}

The process of augmenting a given valuation can of course be iterated, giving rise to the notion of an {\em inductive valuation}. The main result of \cite{MacLane36} is that every MacLane pseudovaluation is inductive.

\begin{defn} \label{def:inductive_valuation}
  A MacLane pseudovaluation $v$ is called {\em inductive} if there exists a chain of MacLane pseudovaluations
  \[
      \vg=v_0, v_1,\ldots,v_n=v,
  \]
  such that $v_i$ is a proper augmentation of $v_{i-1}$, for $i=1,\ldots,n$. We call $v_0,\ldots,v_n$ an {\em augmentation chain} representing $v$.
\end{defn}

We write an inductive valuation as in Definition \ref{prop:inductive_valuations} as
\begin{equation} \label{eq:inductive_valuations2}
v = [v_0,v_1(\phi_1)=\lambda_1,\ldots,v_n(\phi_n)=\lambda_n].
\end{equation}
Here it is understood that $v_0=\vg$ is the Gauss valuation, that $v_n=v$, that $\phi_i$ is a key polynomial for $v_{i-1}$, $\lambda_i> v_{i-1}(\phi_i)$ and that $v_i=[v_{i-1},v_i(\phi_i)=\lambda_i]$. Also, $\lambda_i\neq \infty$ for $i<n$.
We say that the augmentation chain $v_0,\ldots,v_n$ is {\em minimal} if
\begin{equation} \label{eq:minimal_aug_chain}
  \deg\phi_1<\ldots <\deg \phi_n.
\end{equation}

\begin{prop} \label{prop:inductive_valuations}
	Let $v\in V(K[x])^*$. Then there exists a unique minimal augmentation chain representing $v$. In particular, $v$ is inductive.
\end{prop}

\begin{proof}
This is the combination of Theorem 8.1 and Theorem 15.3 of \cite{MacLane36}. See also Theorem 4.31 and Theorem 4.33 of  \cite{RuethThesis}.
\end{proof}

\begin{defn} \label{def:predecessors}
	Given $v\in V(K[x])^*$, let $v_0,\ldots,v_n$ be the unique minimal augmentation chain representing $v$. The valuations $v_0,\ldots,v_{n-1}$ are called the {\em predecessors} of $v$. We set
	\[
	   P(v) := \{v_0,\ldots,v_{n-1}\}.
	\]
	We call $v_{n-1}$ the {\em immediate predecessor} of $v$. Note that $v_0=\vg$ is the Gauss valuation, by definition.
\end{defn}

\begin{lem} \label{lem:predecessors}
	Let $v,w \in V(K[x])^*$ be MacLane pseudovaluations such that $v\leq w$. Then $P(v)\subset P(w)$.
\end{lem}
\begin{proof}
This follows from \cite[Corollary 4.37]{RuethThesis}.
\end{proof}

\subsection{Residue classes and discoids} \label{subsec:residue_classes}

Let $v\in V(K[x])$ be a MacLane valuation. The goal in this subsection is to further analyze the set
\[
   D_v := \{ w\in V(K[x])^* \mid v\leq w\}.
\]
We set $D_v^\circ:=D_v\backslash\{v\}$.

Let $\OO_v \subset K(x)$ be the valuation ring of $v$. We consider the integral domain
\[
A_v := \{ f \in K[x] \mid v(f)\geq 0 \} = K[x]\cap \OO_v.
\]
This is a normal domain in its fraction field $K(x)$.

\begin{lem}
  The prime ideals of height one of $A_v$ are the following:
  \begin{itemize}
  \item
    $\p_v := \{ f \in K[x] \mid v(f)> 0 \}$, and
  \item
    $\p_g:=A_v\cap (K[x]\cdot g)$, where $g\in K[x]$ runs over all irreducible and monic polynomials.
  \end{itemize}
  Moreover, the localization of $A_v$ at any one of these prime ideals is a discrete valuation ring. In particular,
  \[
    (A_v)_{\p_v} = \OO_v.
  \]
\end{lem}

\begin{proof}
Let $v_g$ denote the discrete valuation on $K(x)$ corresponding to an irreducible and monic polynomial $g\in K[x]$. Its valuation ring is the localization of $K[x]$ at the maximal ideal $(g)\lhd K[x]$, $\OO_{v_g}=K[x]_{(g)}$. Since $K[x]$ is a principal ideal domain, we have  $K[x]=\cap_g \OO_{v_g}$. Therefore,
\begin{equation} \label{eq:Krull_intersection}
    A_v = \left(\bigcap_g \OO_{v_g}\right) \cap \OO_v.
\end{equation}
It follows that $A_v$ is a {\em Krull ring}, see \cite[\S 12]{MatsumuraCRT}.

It is easy to see that the family of discrete valuation rings $(\OO_{v_g},\OO_v)$ is a \emph{minimal defining family} for $A_v$, i.e.\ \eqref{eq:Krull_intersection} holds, and none of these rings contains the intersection of the others. By
\cite[Theorem 12.3]{MatsumuraCRT}, $A_v$ has a unique minimal defining family, consisting precisely of the localizations of the prime ideals of $A_v$ of height one. The lemma follows immediately.
\end{proof}

The lemma implies in particular that the field of fractions of the residue ring
\[
\Ab_v := A_v/\p_v
\]
is equal to the residue field $k(v)$ of $v$.
It is proved in \cite[\S 4.1.3]{RuethThesis} that
\begin{equation} \label{eq:Ab_v}
    \Ab_v \cong k'[t]
\end{equation}
is a polynomial ring in one variable over a finite extension $k'/k$.

\begin{lem} \label{lem:residue_class1}
\begin{enumerate}
\item
  For $w\in D_v^\circ$ the set
  \[
      \m_w := A_v\cap\p_w = \{ f\in K[x] \mid v(f)\geq 0,\; w(f)>0 \}
  \]
  is a maximal ideal of $A_v$ containing $\p_v$.
\item
  The map
  \begin{equation} \label{eq:D_v_residue_map}
      w \mapsto \bar{\m}_w :=\m_w/\p_v
  \end{equation}
  is a surjective map from $D_v^\circ$ onto the set of all maximal ideals of $\Ab_v$.
\end{enumerate}
\end{lem}

\begin{proof}
Let $w\in D_v^\circ$, i.e.\ $w>v$. Then obviously $\m_w$ is a prime ideal of $A_v$ properly containing $\p_v$. Therefore, $\bar{\m}_w=\m_w/\p_v$ is a nonzero prime ideal of $\Ab_v$. Since $\bar{A}_v$ is a polynomial ring over a field, $\bar{\m}_w$ is actually a maximal ideal. It follows that $\m_v$ is maximal, and (i) is proved.

Let $\bar{\m}\lhd\Ab_v$ be a maximal ideal and $\m\lhd A_v$ its inverse image. By \eqref{eq:Ab_v}, $\Ab_v$ is a principal ideal domain, so $\bar{\m}=(\bar{f})$, where $\bar{f}$ is the image of a polynomial $f\in\m$. It follows from \cite[Theorem 13.1]{MacLane36} that there exists a key polynomial $\phi$ for $v$ such that $\phi|_v f$. We set
\[
    w := [v,w(\phi)=\infty] \in V(K[x])_\infty.
\]
Then $v<w$ and $0=v(f)<w(f)$. We conclude that
\[
    \m_w=A_v\cap\p_w =\m,
\]
as both sides are maximal ideals which contain $\p_v$ and the element $f$. This proves (ii).
\end{proof}

\begin{defn} \label{def:residue_classes}
  The fibers of the map \eqref{eq:D_v_residue_map} are called the {\em residue classes} of $D_v$. We write $D_v(\bar{\m})$ for the residue class corresponding to a maximal ideal $\bar{\m}\lhd\Ab_v$, and $D_v(w)$ for the residue class containing an element $w\in D_v^\circ$.
\end{defn}

\begin{rem} \label{rem:residue_classes_key_polynomial}
  Let $\phi$ be a key polynomial for $v$.
  Then all proper augmentations
  \[
      v' := [v,v'(\phi)=\lambda],
  \]
  with $\lambda>v(\phi)$, lie in the same residue class, which we denote by $D_v(\phi)$. Using Remark \ref{rem:augmented_valuations} one easily checks that two key polynomials define the same residue class if and only if they are $v$-equivalent. Moreover, it follows from \cite[Theorem 13.1]{MacLane36} that every residue class is of the form $D_v(\phi)$, for some key polynomial $\phi$.
\end{rem}

\begin{lem} \label{lem:residue_class2}
  Let $w_1,w_2\in D_v^\circ$. If $w_1\leq w_2$ then $w_1,w_2$ lie in the same residue class.
\end{lem}

\begin{proof}
If $w_1\leq w_2$ then
\[
   \m_{w_1} = A_v\cap\p_{w_1} \subset A_v\cap\p_{w_2} = \m_{w_2}.
\]
Since both ideals are maximal, we actually have equality.
\end{proof}

\begin{defn} \label{def:discoid}
	Let $g\in K[x]$ be monic, integral and irreducible, and $t\in\QQ$, $t\geq 0$. We set
	\[
	D(g,t):= \{ v\in V(K[x])^* \mid v(g)\geq t \}, \quad
	D(g,t)^\circ:= \{ v\in V(K[x])^* \mid v(g)> t \}.
	\]
	A subset $D\subset V(K[x])^*$ of the form $D(g,t)$ (resp.\ $D^\circ(g,t)$) as above is called a closed (resp.\ an open) {\em discoid}.
\end{defn}

\begin{prop} \label{prop:discoid1}
	Every closed discoid $D$ is of the form
	\[
	D = D_v :=\{ w\in V(K[x])^* \mid v\leq w\},
	\]
	where $v$ is the minimal element of $D$.
\end{prop}

\begin{proof}
This is proved in \cite[\S 4.4]{RuethThesis}.
\end{proof}

\begin{lem} \label{lem:discoid1}
  Let $g\in K[x]$ be monic, integral and irreducible, and let $v_\xi$ be the MacLane pseudovaluation such that $v_\xi(g)=\infty$. Then
  \[
       D_v = D(g,t), \quad \text{where $t:=v(g)$},
  \]
  for any MacLane valuation $v\leq v_\xi$.
  Moreover, the residue class containing $v_\xi$ is
  \[
      D_v(v_\xi) = D^\circ(g,t).
  \]
\end{lem}

\begin{proof}
By Proposition \ref{prop:discoid1} we have
\[
    D(g,t)=D_{v_0},
\]
where $v_0$ is the minimal element of the discoid $D(g,t)$. By construction and the assumption we have $v_0\leq v<v_\xi$, and these are all elements of $D(g,t)$. It follows that $v_0(g)=v(g)=t$.
We can find $a,b\in\NN$ such that $v_0(f)=v(f)=0$, where $f:=g^a/\pi_K^b$. But $v_\xi(f)=\infty >0$, so $f\in\m\backslash\p_v$, where $\m\lhd A_v$ is the maximal ideal corresponding to the residue class $D_{v_0}(v_\xi)$. Since $v(f)=v_0(f)=0$, $v$ does not lie in the residue class $D_{v_0}(v_\xi)$. On the other hand, $v_0\leq v<v_\xi$. Now Lemma \ref{lem:residue_class2} shows that $v=v_0$, proving the first part of the proposition.

By construction, $\p_g:=A_v\cap (K[x]\cdot g)$ is the only minimal prime ideal of $A_v$ containing $f$. Therefore, $\m$, which contains $\p_g$, is the only maximal ideal of $A_v$ containing $f$. This implies
\[
  D_v(v_\xi) = \{ w\in D_{v_0} \mid w(f) >0 \}.
\]
But for all $w\in D_v$ we have $w(f)>0$ if and  only if $w(g)>t$. This proves the second part of the proposition, and we are done.
\end{proof}

\begin{lem} \label{lem:discoid2}
  Let $v\in V(K[x])$ be a MacLane valuation.
  Let $g\in K[x]$ be monic and irreducible (but not necessarily integral) and $\xi:=(g)\in|\PP_K^1|$. Let $a,b\in\ZZ$, with $a>0$, such that
  \[
     v(f) = 0, \quad \text{where}\;\; f:=\frac{g^a}{\pi_K^b}.
  \]
  Let $\bar{f}\in \bar{A}_v$ be the image of $f$. Then the following are equivalent:
  \begin{enumerate}[(a)]
  \item
    $\bar{f}\in\Ab_v$ is not a unit,
  \item
    $v<v_\xi$.
  \end{enumerate}
\end{lem}

\begin{proof}
The implication (b)$\Rightarrow$(a) has been shown during the proof of Proposition \ref{lem:discoid1}. To prove (a)$\Rightarrow$(b), we assume that $\bar{f}$ is not a unit. As in the proof of Lemma \ref{lem:residue_class1} we conclude that there exists a MacLane pseudovaluation $w$ such that $v\leq w$ and $v(f)<w(f)$. Then we also have $v(g)<w(g)$.

Set $t:=v(g)$ and let $v_0$ be the minimal element of the discoid $D(g,t)$, such that $D_{v_0}=D(g,t)$ (Lemma \ref{lem:discoid1}). Since $w(g)>v(g)=t$, Lemma \ref{lem:discoid1} also shows that $w$ and $v_\xi$ lie in the same residue class. If $v=v_0$, then clearly $v\leq v_\xi$ and we are done. Otherwise, $v<w$ and Lemma \ref{lem:residue_class2} shows that $v$, $w$ and $v_\xi$ all lie in the same residue class of $D_{v_0}$. But then $v(g)>t$ by Lemma \ref{lem:discoid1}, which is a contradiction. This completes the proof of the lemma.
\end{proof}

\begin{lem} \label{lem:residue_classes3}
	Two elements $w_1,w_2\in D_v^\circ$ lie in the same residue class if and only if there exists $u\in D_v^\circ$ such that $u\leq w_1,w_2$.
\end{lem}

\begin{proof}
Suppose that there exists $u\in D_v^\circ$ with $u\leq w_1,w_2$. Then Lemma \ref{lem:residue_class2} immediately implies that $w_1,w_2$ lie in the same residue class. Conversely, assume that $w_1,w_2$ both lie in the same residue class as the infinite pseudovaluation $v_\xi$, with $\xi=(g)$. Lemma \ref{lem:discoid1} implies that
\[
   t:=v(g) < t_1:=\min(w_1(g),w_2(g)).
\]
Let $u$ be the minimal element of the discoid $D(g,t_1)$.
Then $D_u=D(g,t_1)$, by Proposition \ref{prop:discoid1}. Therefore, $v<u\leq w_1,w_2$. This completes the proof of the lemma.
\end{proof}

\subsection{Valuation trees} \label{subsec:valuation_tree}

Using the results of the previous subsection, we will now prove that any finite subset of $V(K[x])^*$ can be extended in such a way that it forms a finite rooted tree.

\begin{prop}
  The subset
  \[
     \{w \mid w\leq v \} \subset V(K[x])^*
  \]
  is totally ordered, for all $v\in V(K[x])^*$.
\end{prop}

\begin{proof}
It suffices to prove this when $v=v_\xi$ is an infinite pseudovaluation. If $w\leq v_\xi$ then
\[
   D_w = D(g,t),
\]
where $g$ is the irreducible polynomial corresponding to $\xi$ and $t:=w(g)$ (Lemma \ref{lem:discoid1}). Given another valuation $w'\leq v_\xi$, we conclude that
\[
    w\leq w' \quad\Leftrightarrow\quad w(g)\leq w'(g).
\]
Therefore, any two valuations $w,w'\leq v_\xi$ are comparable.
\end{proof}

\begin{prop}
  For any two elements $v,w\in V(K[x])^*$, the set
  \[
     \{ u\in V(K[x])^* \mid u\leq v,\; u\leq w \}
  \]
  has a unique maximal element.
\end{prop}

\begin{proof}
Choose an infinite pseudovaluation $v_\xi$ such that $v\leq v_\xi$, and let $g$ denote the corres\-ponding irreducible polynomial. Set
\[
    t := \min\{ v(g), w(g) \}
\]
and let $u$ be the minimal element of the discoid $D(g,t)$. Then
\[
    D_u = D(g,t),
\]
by Lemma \ref{lem:discoid1}. In particular, $u\leq v,w$.

Assume that there exists $u'\in V(K[x])^*$ such that $u<u'\leq v,w$. Then Lemma \ref{lem:residue_class2} shows that $u',v,w$ all lie in the same residue class of the discoid $D_u$ as $v_\xi$. But then Lemma \ref{lem:discoid1} shows that $v(g),w(t)>t$, which contradicts the definition of $t$. We conclude that $u$ is the desired maximal element. The uniqueness  is clear.
\end{proof}

\begin{defn} \label{def:inf}
  Given $v,w\in V(K[x])^*$, the maximal elements $u$ from the previous lemma is called the {\em infimum} of $v,w$,
  \[
     \inf(v,w) := \max \{ u \mid u\leq v,w \}.
  \]
\end{defn}

\begin{cor}
  Let $V\subset V(K[x])^*$ be a finite nonempty subset. Assume that $V$ is {\em inf-closed}, i.e.\ for $v,w\in V$ we also have $\inf(v,w)\in V$. Then $V$, as a partially ordered set, is a {\em rooted tree}, i.e.\ for all $v\in V$ the subset
  \[
      \{w\in V \mid w\leq v\}
  \]
  is well ordered, and there exists a minimal element $v_0\in V$ (the {\em root}).
\end{cor}

\begin{rem} \label{rem:valuation_tree}
 A finite, inf-closed subset $V\subset V(K[x])^*$ is called a {\em valuation tree}. We may consider a valuation tree as a finite rooted tree in the sense of graph theory, as follows. The set of (directed) edges is the set of pairs $(v,w)$, with $v,w\in V$, such that $v<w$ and there is no element $u\in V$ with $v< u< w$. The two valuations $v,w$ (nodes of $V$) are then said to be {\em adjacent}. The root of $V$ is the minimal element $v_0$.
\end{rem}

\subsection{Algorithmic tools} \label{subsec:algorithmic_tools}

The results presented in this section have been implemented in the computer algebra system {\tt Sage} (\cite{sagemath}) and the Sage toolbox {\tt MCLF} (\cite{mclf}). The algorithms described in \S \ref{sec:regular_models_of_projective_line} and \S \ref{sec:global_sections} crucially rely on these implementations.

{\tt Sage} contains a very general concept of discrete (pseudo)valuations on rings and fields. See the relevant chapter of the \href{http://doc.sagemath.org/html/en/reference/valuations/index.html}{documentation}. For instance, it is possible to create and compute with inductive valuations $v$ on a polynomial ring $K[x]$ and the function field $K(x)$, where $K$ is an arbitrary field equipped with a discrete (and possibly trivial) valuation $v_K$. Given $v$ it is possible to construct an arbitrary augmentation
\[
    v' = [v, v'(\phi)=\lambda]
\]
(\S \ref{subsec:augmentations}). One can evaluate $v(f)$, for polynomials and rational functions $f$, and compute the reduction $\bar{f}\in k(v)$ if $v(f)\geq 0$. Given an irreducible element of the residue ring $\psi\in \Ab_v=k'[t]$, one can `lift` $\psi$ to a key polynomial $\phi$ for $v$ (in the situation of Remark \ref{rem:residue_classes_key_polynomial}, find $\phi$ such that $D_v(\phi)=D_v(\m)$, where $\bar{\m}:=\m/\p_v=(\psi)$). Finally, given a square free polynomial $f\in K[x]$, it is possible to find arbitrarily good approximations $v$ for the MacLane pseudovaluations $v_\xi$ corresponding to the irreducible factors $g\mid f$, $g\in\hat{K}[x]$, over the completion $\hat{K}$ of $K$ with respect to $v_K$. When $f=g$ is already irreducible over $\hat{K}$, then the algorithm computes the inductive pseudovaluation
\[
    v_\xi = [v_0,v_1(\phi_1)=\lambda_1,\ldots, v_n(g)=\infty].
\]
In the general case, the result of the algorithm can be used to find an arbitrarily good approximate factorization of $f$ over $\hat{K}$.

\begin{rem} \label{rem:completeness_of_K}
  So while we always assume in this article that $K$ is complete with respect to $v_K$, we do not actually work with complete fields at all when doing explicit computations. This is not a serious restriction for applications, and it makes it much easier to obtain results which are provably correct.
\end{rem}

The results of \S \ref{subsec:residue_classes} on discoids and residue classes, and of \S \ref{subsec:valuation_tree} on valuation trees have been implemented within the module \href{https://mclf.readthedocs.io/en/latest/berkovich.html}{ \tt berkovich}, which is a part of the Sage toolbox {\rm MCLF}. The main object one can work with is the {\em Berkovich line} over a discretely valued field $(K,v_K)$; it represents the analytification $X^\an$ of the projective line $X=\PP^1_{\hat{K}}$ over the completion $\hat{K}$ of $K$, see Remark \ref{rem:berkovich}. There are three types of points on this space one can use: points of Type I (corresponding to infinite pseudovaluations on $F_X:=K(x)$), Type II (corresponding to true valuations on $F_X$) and Type V. The latter are not actual points on the analytic space $X^\an$, but rather on the {\em adic space} $X^{\rm ad}$. They correspond to certain valuations of rank $2$ and, in the language of \S \ref{subsec:residue_classes} to residue classes of discoids. These points form a partially ordered set (of which $V(K[x])^*$ is a subset), corresponding to the closed unit disk. Making essential use of the results of \S \ref{subsec:residue_classes} and \S \ref{subsec:valuation_tree}, it is possible to evaluate inequalities, compute infima (Definition \ref{def:inf}) and build the valuation tree spanned by a given finite set of points (Remark \ref{rem:valuation_tree}).

%% file: regular_models_of_P1.tex

\section{Regular models of the projective line} \label{sec:regular_models_of_projective_line}

We fix a discretely valued field $(K,v_K)$ as in the previous section and a smooth projective curve $X_K$ over $K$. Let $F_X$ denote the function field of $X_K$. By a {\em model} of $X_K$ we will mean a normal $\oo_K$-model, i.e.\ an $\oo_K$-scheme $X\to\Spec\,\oo_K$ which is proper and flat and whose generic fiber is equal to $X_K$. In this section we will mostly be interested in the case $X_K=\PP^1_K$. We then write $F_X=K(x)$.

Let $D_K\subset X_K$ be an effective, reduced divisor (in other words, $D_K$ is a finite set of closed points of $X_K$). For any model $X$ of $X_K$, the closure $D^\hor\subset X$ of $D_K$ is a closed subscheme of dimension $1$, flat over $\Spec\,\oo_K$. We call $D^\hor$ a {\em horizontal divisor} on the model $X$.

The following result is well known (combine, for instance, the main result of \cite{Lipman78} with \cite[Theorem 9.2.26.]{liu2002algebraic}).

\begin{prop} \label{prop:regular_model_of_P1}
  Let $D_K\subset X_K$ be an effective reduced divisor. There exists a model $X$ of $X_K$ with the following properties.
  \begin{enumerate}
  \item
    $X$ is regular.
  \item
    Let $D^\hor\subset X$ be the closure of $D_K$ inside $X$ and $D^{\rm vert}:=X_s\subset X$ the reduced special fiber. Then $D:=D^{\rm vert}\cup D^\hor$ is a normal crossing divisor on $X$.
  \end{enumerate}
\end{prop}

The goal of this section is to make the proof of Proposition \ref{prop:regular_model_of_P1} explicit in the case where $X_K=\PP^1_K$, and describe an algorithm for computing such a model $X$ (see \S \ref{subsec:algorithm}). We also describe an actual implementation of this algorithm.
We use the methods of \cite{RuethThesis} and  \cite{ObusWewers}. See \cite{ObusSrinivasan} for a slightly different treatment of essentially the same results.

\subsection{Models and valuations}

Given a normal model $X$ of $X_K$, we let $X_s$ denote the special fiber. Here we consider $X_s$ as a closed subscheme of $X$ with the reduced subscheme structure. We can write $X_s$ as a union of irreducible components, which we usually call $E_i$,
\[
     X_s = \cup_i E_i.
\]
So each $E_i\subset X$ is a prime Weil divisor and thus gives rise to a discrete valuation
\[
    v_i=v_{E_i} :F_X \to\QQ\cup\{\infty\}.
\]
The restriction of $v_i$ to the base field $K$ is equivalent to $v_K$. We normalize $v_i$ in such a way that we have equality, $v_i|_K=v_K$. Then the value group of $v_i$ is of the form $v_i(F_X^\times)=m_i^{-1}\ZZ$ for a positive integer $m_i$ (recall that $v_K(K^\times)=\ZZ$). We call $m_i$ the {\em multiplicity} of the component $E_i$ in the special fiber.

A discrete valuation $v$ on $F_X$ is called {\em geometric} if $v|_K=v_K$ and the residue field $k(v)$ has transcendence degree one over $k$. Let $V(F_X)$ denote the set of all geometric valuations.
Given a normal $\oo_K$-model $X$ of $X_K$, we write $V(X)$ for the finite nonempty set of geometric valuations corresponding to the irreducible components of $X_s$. Then the association
\[
    X \mapsto V(X)
\]
is an order-reversing bijection between the set of all models of $X_K$ and the set of all finite, nonempty subsets of $V(F)$ (\cite[Corollary 3.18]{RuethThesis}). This will allow us to construct models of $X$ in an easy way.

It will be useful to consider the following enhancement of the above bijection. Let $\abs{X_K}$ denote the set of all closed points of $X_K$. For $\xi\in \abs{X_K}$ we let $\OO_\xi\subset F_X$ denote the local ring of $\xi$. We define a map
\[
   v_\xi:\; \left\{ \begin{array}{ccc}
     F_X  & \longrightarrow & \QQ\cup\{\pm\infty\} \\
     f = g/h & \mapsto & v_K(g(\xi)) - v_K(h(\xi)).
   \end{array} \right.
\]
Here $g,h\in\OO_\xi$ are chosen such that not both lie in the maximal ideal of $\OO_\xi$. We use $g(\xi), h(\xi)$ as a  suggestive notation for the images of $g,h$ in the residue field $k(\xi)$ of the valuation ring $\OO_\xi$ (which is a finite extension of $K$), and $v_K$ is used here for the unique extension of $v_K$ to $k(\xi)$.
We call such a map $v_\xi$ a {\em geometric infinite pseudovaluation} on $F_X$.

We note that $v_\xi|_K = v_K$, and that $v_\xi$ uniquely determines the point $\xi\in X_K$. For instance,
\[
    \OO_\xi = \{ f\in F_X \mid v_\xi(f)> -\infty \},\quad
    \m_\xi = \{ f\in F_X \mid v_\xi(f)= \infty \}.
\]
But $v_\xi$ should not be confused with the discrete valuation corresponding to the valuation ring $\OO_\xi$.

We set
\[
    V(F_X)^* := V(F_X)\cup \{v_\xi \mid \xi\in\abs{X_K} \}.
\]
The elements of this set are called {\em geometric pseudovaluations}. The subset of infinite pseudovaluations is denoted $V(F_X)_\infty$.

Given a pair $(X,D^\hor)$ consisting of a model $X$ of $X_K$ and a horizontal divisor $D^\hor\subset X$, we set
\[
   V^*(X,D):= V(X) \cup \{ v_\xi \mid \xi \in D_K\}.
\]
We have the following formal consequence of \cite[Corollary 3.18]{RuethThesis}.

\begin{prop} \label{prop:models_valuations}
  The association
  \[
      (X,D^\hor)  \mapsto V(X,D^\hor)
  \]
  defines a bijection between the following two sets:
  \begin{itemize}
  \item
    Pairs $(X,D^\hor)$, where $X$ is a model of $X_K$ and $D^\hor\subset X$ a horizontal divisor.
  \item
    Finite subsets of $V(F_X)^*$ with nonempty intersection with $V(F_X)$.
  \end{itemize}
\end{prop}

From now on, we set $X_K:=\PP^1_K$ and hence $F_X=K(x)$ is the rational function field in $x$. Then there is a certain subset of  $V(F_X)^*$  corresponding to MacLane pseudovaluations, as defined in the previous section. More precisely:

\begin{rem} \label{rem:maclane_valuations_extend}
\begin{enumerate}
\item
  Let $v\in V(K[x])^*$ be a MacLane pseudovaluation (Definition \ref{def:maclane_valuation}). Then
  \[
      v\left(\frac{f}{g}\right)= v(f) - v(g), \qquad \text{with $f,g\in K[x]$, $\gcd(f,g)=1$,}
  \]
  defines an extension of $v$ to a geometric pseudovaluation $v:F_X\to\QQ\cup \{\pm\infty\}$. We obtain an inclusion
  \[
      V(K[x])^* \inj V(F_X)^*
  \]
  whose image is the subset of geometric pseudovaluations $v$ such that $v(x)\geq 0$. We will henceforth consider $V(K[x])^*$ as a subset of $V(F_X)^*$.
\item
  It is clear that
  \[
	V(K(x))^* = V(K[x])^* \cup V(K[x^{-1}])^*,
  \]
  so every geometric pseudovaluation is a MacLane pseudovaluation either with respect to $x$ or to $x^{-1}$. As we will see, restricting attention to the subset $V(K[x])^*$ is not a serious restriction, and it considerably simplifies the exposition. In our implementation (\cite{regular_models}) this restriction is avoided.
\end{enumerate}
\end{rem}

\subsection{Explicit description of the divisor $D$}
\label{subsec:component_tree}

Let us fix a finite set $V^*\subset V(F_X)^*$ containing at least one finite valuation. Let $(X,D^\hor)$ be the model of $(X_K=\PP^1_K, D_K)$ corresponding to $V^*$ via Proposition \ref{prop:models_valuations}. Let $X_s$ denote the special fiber of $X$, considered as a closed subscheme with its reduced subscheme structure. Finally, set $D:=X_s\cup D^\hor$.

By construction, the irreducible components of $D$ correspond to the elements of $V^*$,
\[
   D = \bigcup_{v\in V^*} E_v.
\]
Here $E_v\subset X_s$ for $v\in V:=V^*\cap V(F_X)$ and $E_v\subset D^\hor$ for $v\in V_\infty:=V^*\cap V(F_X)_\infty$. The goal of this subsection is to describe the scheme $D$ explicitly, at least as a topological space, using only the set $V^*$. We make the following assumptions on $V^*$.

\begin{ass} \label{ass:V}
\begin{enumerate}[(a)]
\item
  All elements of $V^*$ are MacLane pseudovaluations, $V^*\subset V(K[x])^*$. Thus, we may consider $V^*$ as a partially ordered set (see Definition \ref{def:maclane_valuation} for the definition of $\leq$).
\item
  For $v_1,v_2\in V^*$, $\inf(v_1,v_2)$ is also contained in $V^*$ (Definition \ref{def:inf}).
\item
    Given $v\in V^*$, every predecessor of $v$ is contained in $V$ as well (i.e.\ $P(v)\subset V$, see Definition \ref{def:predecessors}).
\end{enumerate}
\end{ass}

Assumptions (a) and (b) mean that we may regard $V^*$ as a rooted tree, see Remark \ref{rem:valuation_tree}.
Recall that a pair $(v_1,v_2)$ of elements of $V^*$ is an oriented edge of this tree if $v_1<v_2$ and if there is no element $v\in V$ such that $v_1<v<v_2$. If this is the case, we call the two elements $v_1,v_2\in V^*$ {\em adjacent}.

The main result of this subsection states that the tree $V^*$ is naturally isomorphic to the {\em tree of components} given by the divisor $D$. Assumption \ref{ass:V} (c) is not needed for this result and will only be used later in \S \ref{subsec:regularity}.

\begin{prop} \label{prop:component_tree}
  Let $V^*$, $X$, $D$, $E_v$ be as before. Assume that Assumption \ref{ass:V} (a) and (b) hold.
  Then for all $v,v'\in V^*$, $v\neq v'$, the components $E_v$ and $E_{v'}$ intersect if and only if $v,v'$ are adjacent. If this is the case then there is a unique intersection point.
\end{prop}

\begin{proof}
The proof will occupy most of the rest of this subsection. We start with two lemmata.

\begin{lem} \label{lem:component_tree1}
  Let $v\in V(K[x])$ be a MacLane valuation, and let $X_v$ denote the model of $X_K$ corresponding to the set of valuations $\{v\}$. Let $\overline{\infty}\subset X_v$ denote the closure of the point $\infty\in X_K=\PP^1_K$, and $U_v:=X_v \backslash \{\overline{\infty}\}$.
  Write $\infty_v\in X_{v,s}$ for the unique intersection point of $\overline{\infty}$ with $X_{v,s}$.
  \begin{enumerate}
  \item
    $U_v=\Spec\, A_v$, where $A_v:= \OO_v\cap K[x]$.
  \item
    $X_{v,s}-\{\infty_v\}=\Spec \bar{A}_v$, where
    $\bar{A}_v:=A_v/\p_v$, and where $\p_v:=\{f\in K[x]\mid v(f)>0\}$.
  \item
    Let $v_\xi\in V(F_X)_\infty$ be an infinite geometric pseudovaluation, corresponding to a closed point $\xi\in X_K$. Let $\bar{\xi}\subset X_v$ denote the closure of $\xi$. Then $\bar{\xi}$ intersects $X_{v,s}$ in $\infty_v$ if $v_\xi(x) <0$, or if $v_\xi\in V(K[x])_\infty$ but $v\not<v_\xi$. Otherwise, $\bar{\xi}$ intersects $X_{v,s}$ in the point $\m\in\Spec \bar{A}_v$ corresponding to the residue class of the discoid $D_v$ containing $v_\xi$ (Definition \ref{def:residue_classes}).
  \end{enumerate}
\end{lem}

\begin{proof}
By construction, $X_v\to\Spec\,\oo_K$ is a projective fibred surface whose reduced special fiber $X_{v,s}$ is a projective irreducible curve over $k$. The composition $\overline{\infty}\inj X_v\to\Spec\, \oo_K$ is finite and birational. Since $\oo_K$ is integrally closed, it follows that $\overline{\infty}\iso\Spec\,\oo_K$ is actually an isomorphism. Using \cite[Corollary 5.3.24]{liu2002algebraic} we see that the divisor $\overline{\infty}$ is ample on $X_v$. Therefore, the complement $U_v=X_v\backslash \overline{\infty}$ is an affine open subset. So $U_v=\Spec A_v$, where $A_v\subset F_X=K(x)$ is the subring of functions which are regular on $U_v$. By definition of $X_v$, these are precisely the functions $f\in K(x)$ such that
\[
    v(f)\geq 0, \quad\text{and}\quad \ord_\xi(f)\geq 0\;\;\forall \xi\in|X_K|,\,\xi\neq \infty.
\]
The second condition is equivalent to $f$ being a polynomial in $x$. This proves (i). Part (ii) is an immediate consequence of (i) and the definitions.

Let $\xi\in|X_K|$ be a closed point on the generic fiber, $\xi\neq\infty$. Then $\xi$ corresponds to a maximal ideal of $K[x]$, which is generated by a monic irreducible polynomial $g$. Set
\[
    f :=\frac{g^a}{\pi^b}, \quad \text{where $v(g)=b/a$, $a>0$.}
\]
Then $v(f)=0$. In particular, $f\in A_v$. We let $\bar{f}\in\bar{A}_v$ denote the image of $f$. By construction, the principal divisor of $f$ as a rational function on $X_v$ is
\[
   (f) = a\cdot\bar{\xi} - \deg(f)\cdot\overline{\infty}.
\]
As $\X$ is proper over $\Spec\OO_K$ and $\OO_K$ is complete, $\bar{\xi}$ intersects the special fiber $X_{v,s}$ in a unique closed point (the {\em reduction} of $\xi$, see \cite[\S 10.1.3]{liu2002algebraic}).
Now it follows from (ii) that this point of intersection is $\infty_v$ if and only if $\bar{f}\in\bar{A}_v$ is a unit. Otherwise, the point of intersection corresponds to the maximal ideal
\[
    \m :=\rad(\bar{f})\lhd\bar{A}_v.
\]
But Lemma \ref{lem:discoid2} says that $\bar{f}$ is not a unit if and only if $v<v_\xi$. Moreover, if this is the case then the ideal $\m$ defined above corresponds precisely to the residue class $D_v^\circ(\xi)$, as in Definition \ref{def:residue_classes}. Statement (iii) of the lemma follows.
\end{proof}

\begin{lem} \label{lem:component_tree2}
  Fix one element $v\in V$.  Let
  \[
     \phi_v:X\to X_v
  \]
  be the morphism of models corresponding to the inclusion $\{v\}\subset V$. Then the following holds.
  \begin{enumerate}
  \item
    The restriction of $\phi_v$ to $E_v\subset  X$ is a homeomorphism $E_v\iso X_{v,s}$.
  \item
    Let $v'\in V\backslash\{v\}$. Then $\phi_v$ contracts the vertical component $E_{v'}$ to a closed point $z$ of $X_{v,s}$.
    We have $z\neq\infty_v$ if and only if $v<v'$. Moreover, if this is the case then $z$ is the point on $X_{v,s}=\Spec\,\bar{A}_v$ corresponding to the residue class $D_v^\circ(v')$.
  \item
    Let $v'\in V_\infty$. Then the horizontal divisor $\phi(E_{v'})\subset X_v$ intersects the special fiber $X_{v,s}$ in a point $z\neq\infty_v$ if and only if $v<v'$. Moreover,  if this is the case then $z$ is the point on $X_{v,s}=\Spec\,\bar{A}_v$ corresponding to the residue class $D_v^\circ(v')$.
  \end{enumerate}
\end{lem}

\begin{proof}
By construction, the morphism $\phi_v:X\to X_v$ contracts all vertical components $E_{v'}$, where  $v'\neq v$, to closed points of $X_{v,s}$, and it induces a finite, birational map $E_v\to X_{v,s}$. Both $E_v$ and $X_{v,s}$ are proper and integral curves over $k$. By Lemma \ref{lem:component_tree1} (ii), the open subset $X_{v,s}\backslash\{\infty_v\}$ is isomorphic to the affine line over a finite extension $k'/k$. It follows that the smooth projective model of both $E_v$ and $X_{v,s}$ is the projective line $\PP^1_{k'}$ and that the natural birational morphism $\PP^1_{k'}\to X_{v,s}$ is a homeomorphism. since this map factors through the curve $E_v$, the map $E_v\to X_{v,s}$ is a homeomorphism, too, which proves (i).

A vertical component $E_{v'}$, with $v'\neq v$, is contracted to the center of the discrete valuation $v'$ on the scheme $X_v$.  Clearly, this center is a closed point on $X_{v,s}$. It lies on the affine open $\Spec A_v$ if and only if
\[
    A_v\subset\OO_{v'}.
\]
But this latter condition holds if and only if $v< v'$, by definition of the order relation $\leq $. If this is the case, then the center of $v'$ on $\Spec A_v$ is the maximal ideal
\[
    \m :=A_v\cap \p_{v'}
\]
which, by Lemma \ref{lem:residue_class1}, corresponds to the residue class $D_v(v')$. This proves (ii).
Part (iii) follows directly from Lemma \ref{lem:component_tree1} (iii).
\end{proof}

Now we prove Proposition \ref{prop:component_tree}.
Let $v,v'\in V^*$ be distinct. We have to prove that $E_v$ intersects $E_{v'}$ if and only if $v$ and $v'$ are adjacent, and that there is at most one intersection point. For simplicity, we assume that $v,v'\in V$ are both MacLane valuations, and hence $E_v$ and $E_{v'}$ are vertical components. The argument for the general case is very similar and left to the reader (essentially, one uses Part (iii) of Lemma \ref{lem:component_tree2} instead of Part (ii)).

We first assume that $v,v'$ are {\em not} adjacent. There are two cases to consider. Firstly, suppose that there exists $v''\in V^*$ such that $v<v''<v'$. By Lemma \ref{lem:component_tree2} (ii), the morphism $\phi_{v''}:X\to X_{v''}$ contracts $E_v$ to the point $\infty_{v''}\in X_{v'',s}$ and $E_{v'}$ to a closed point distinct from $\infty_{v''}$. Therefore, $E_v$ and $E_{v'}$ do not intersect. In the second case we assume that $v$ and $v'$ are incomparable (i.e.~neither $v<v'$ nor $v'<v$). By Assumption \ref{ass:V} (b), $v'':=\inf(v,v')\in V$. Moreover, $v,v'$ lie in different residue classes of the discoid $D_{v''}$, by Lemma \ref{lem:residue_classes3}. Hence, by Lemma \ref{lem:component_tree2} (ii), the morphism $\phi_{v''}$ contracts $E_v$ and $E_{v'}$ to distinct points. Therefore, $E_v$ and $E_{v'}$ do not intersect, as in the first case.

Conversely, assume that $E_v$ and $E_{v'}$ do not intersect. We want to show that $v,v'$ are not adjacent. So we may assume that $v<v'$. Since $X_s$ is connected, there exists a chain of pairwise distinct components
\[
   E_v = E_{v_0}, \dots, E_{v_n} = E_{v'}
\]
with $E_{v_i} \cap E_{v_{i-1}} \neq \emptyset$ and $n\geq 2$. By Lemma \ref{lem:component_tree2} (ii), the morphism $\phi_{v}$ contracts each component $E_{v_1},\ldots, E_{v_n}$ to a single point. But since $E_{v_i}\cap E_{v_{i-1}}\neq\emptyset$, they are all contracted to the same point. Since $v<v'$, applying Lemma \ref{lem:component_tree2} again shows that $v<v_1,\ldots,v_n=v'$. Similarly, one can see that $v_1, \dots, v_{n-1}<v'$. In this case, consider the map $\phi_{v'}: X \to X_{v'}$. All components $E_{v_0}, \dots, E_{v_{n-1}}$ get contracted to the same point. Since $v<v'$, this is a point $z \neq \infty_v$, and applying Lemma \ref{lem:component_tree2} again, it follows that $v_1,\dots, v_{n-1}< v'$.

The same reasoning shows that $v_i < v_j$ for all $i<j$ and so
\[
   v=v_0 < \dots < v_n=v'.
\]
Since $n\geq 2$, $v$ and $v'$ are not adjacent.

To finish the proof, we assume that $v<v'$ are adjacent. As we have shown above, $E_v\cap E_{v'}$ is nonempty. However, by Lemma \ref{lem:component_tree2} (i)-(ii), the morphism $\phi_v$ is injective on $E_v$ and contracts $E_{v'}$ to a single point. Therefore, $E_v$ and $E_{v'}$ intersect in a unique point. This concludes the proof of Proposition \ref{prop:component_tree}.
\end{proof}

\subsection{Regularity criteria} \label{subsec:regularity}

We continue with the assumptions and notation of the previous subsection. We will now formulate sufficient conditions for the model $X$ to be regular and $D$ to be a normal crossing divisor. These conditions will be formulated in Lemma \ref{lem:regularity_cond1} and Lemma \ref{lem:regularity_cond2} below.

The claim that $X$ is regular and $D$ is a normal crossing divisor needs only to be checked in a Zariski neighborhood in $X$ of every closed point on $D$. By Proposition \ref{prop:component_tree}, a closed point $z\in D$ lies on one and at most on two vertical components $E_v$, $v\in V$. We write $X_z$ for the localization of $X$ at $z$, $D_z$ for the restriction of $D$ to $X_z$ and $D_z^{\rm vert}$ (resp.\ $D_z^{\rm hor}$) for the vertical (resp.\ the horizontal) part of $D_z$.

\vspace{2ex}\noindent
Let us first consider the case where $z$ lies on exactly one vertical component $E_v$, $v\in V$. Consider the contraction map
\[
     \phi_v:X\to X_v
\]
from Lemma \ref{lem:component_tree2}. Our assumption that $z$ lies on no other vertical component except $E_v$ implies that $\phi_v$ is an isomorphism in a neighborhood of $z$.

Assume first that $\phi_v(z) = \infty_v$. We claim that this implies that $v=\vg$ is the Gauss valuation. To prove this, note that Assumption \ref{ass:V} (b) and (c) imply that the Gauss valuation $\vg$ is an element of $V$. It is in fact the minimal element of $V$. So if $v\neq \vg$ then Lemma \ref{lem:component_tree2} (ii) would imply that $z$ also lies on the component $E_{v'}$, where $v'<v$ is the unique element of $V$ which is adjacent to $v$ and precedes it. This contradicts our assumption, and shows that $v=\vg$ if $\phi_v(z)=\infty_v$.

By the definition of the Gauss valuation, the model $X_{\vg}=\PP^1_{\oo_K}$ is simply the projective line over $\oo_K$, which is smooth over $\oo_K$. It follows that $X_z$ is regular and that $D_z^{\rm vert}$ is a regular, principal divisor, defined by the uniformizer $\pi_K$. Also, $D_z^{\rm hor}=\emptyset$, and $D_z$  is a normal crossing divisor.

So for the rest of this subsection, we may assume that $\phi_v(z)\neq \infty_v$. We use the contraction map $\phi_v:X\to X_v$ to identify $X_z$ with $\Spec\,(A_v)_\m$, where $A_v$ is defined as in \S \ref{subsec:component_tree} and $\m\in\Spec\,A_v\subset X_v$ is the image of $z$. Recall that we defined a {\em residue class} $D_{v}(\m)\subset V(K[x])^*$ corresponding to $\m$, see Definition \ref{def:residue_classes}.

We write $v$ as an inductive valuation,
\begin{equation} \label{eq:regularity_cond1}
   v = [v_0,v_1(\phi_1)=\lambda_1,\dots,v_n(\phi_n)=\lambda_n],
\end{equation}
where $v_0=\vg$, see \eqref{eq:inductive_valuations2}.
 As in Remark \ref{rem:residue_classes_key_polynomial}, we denote by $D_v(\phi_n)$ the residue class associated to $\phi_n$ (which contains e.g.\ the infinite pseudovaluation $v_\xi$, where $v_\xi(\phi_n)=\infty$).

\begin{lem} \label{lem:regularity_cond1}
  Let $z\in X_s$ be a closed point, which lies on exactly one vertical component $E_v$.
  \begin{enumerate}
  \item
    The point $z\in X$ is a regular point of $X$ if and only if one of the following conditions holds:
    \begin{enumerate}[(a)]
    \item
      $n=0$,
    \item
      $\lambda_n\in \gen{1,\ldots,\lambda_{n-1}}_\ZZ$,
    \item
      $D_v(\m)\neq D_v(\phi_n)$.
    \end{enumerate}
  \item
    If $z\in X$ is regular, then $D$ is a normal crossing divisor in a neighborhood of $z$.
  \end{enumerate}
\end{lem}

\begin{proof}
Part (i) is \cite[Lemma 7.3]{ObusWewers}. We recall the proof of the sufficiency of either (a), (b) or (c). Firstly, if (a) holds then $v=\vg$. We have already seen that $X\to \Spec\,\oo_K$ is then smooth in a neighborhood of $z$, and hence $z$ is a regular point.

By construction, the value group of the inductive valuation \eqref{eq:regularity_cond1} is
\[
    v(F_X^\times) = \gen{1,\lambda_1,\ldots,\lambda_n}_\ZZ
      = \frac{1}{e_v}\ZZ
\]
(\cite[Theorem 6.6]{MacLane36}).
This means that there exists a uniformizer for $v$ of the form
\[
    u_v := c\cdot\phi_1^{a_1}\cdot\ldots\cdot \phi_n^{a_n},
    \qquad \text{with}\quad v(u_v) =\frac{1}{e_v},
\]
where $c\in K^\times$ and $a_i\leq 0$. By construction, the principal divisor of $u_v$ on $\Spec A_v$ is
\[
  (u_v) = \p_v + \sum_{i=1}^n a_i\cdot \p_{\phi_i}.
\]
We are interested in the restriction of this divisor to the local scheme $X_z=\Spec (A_v)_\m$.

For $i\neq n$, the prime ideal $\p_{\phi_i}:=A_v\cap (K[x]\cdot\phi_i)\lhd A_v$ is a maximal ideal not containing $\p_v$, by Lemma \ref{lem:discoid2}. For $i=n$, the prime ideal $\p_{\phi_n}$ is contained in $\m$ if and only if $D_v(\m)=D_v(\phi)$, i.e.\ if and only if Condition (c) is false.

In particular, if Condition (c) holds then $D_z$ is a principal divisor, defined by $u_v$. Since $D_z$ is regular, it follows that $X_z$ is regular, too. Similarly, if Condition (b) holds then we may assume that $a_n=0$, and the same conclusion holds. So in both cases, $X_z$ is regular.

It remains to prove (ii).
If $D^{\rm hor}=\emptyset$ then we have already proved that $D_z=D_z^{\rm vert}$ is a regular, and hence a normal crossing divisor. Therefore, we may assume that $z$ lies on a horizontal component $E_{v_\xi}=\bar{\xi}$. Here $\bar{\xi}$ is the closure of a closed point $\xi\in X_K$ corresponding to an infinite pseudovaluation $v_\xi$ which lies in $V^*$. The point $\xi$ corresponds to a monic, integral and irreducible polynomial $g\in K[x]$. Assumption \ref{ass:V} (c) guarantees that the immediate predecessor of $v_\xi$ lies in $V^*$ as well; by Proposition \ref{prop:component_tree}, this predecessor must be the valuation $v$. It follows that $g$ is a key polynomial for $v$ and that
\[
    v_\xi = [v,v_\xi(g)=\infty].
\]
We set
\[
    f := u_v^{-e_v v(g)}\cdot g \in K[x].
\]
Then $v(f)=0$, and hence
\[
    D_z^{\rm hor} = (f)
\]
is a principal divisor. To prove that $D_z$ is normal crossing, we need to show that $\bar{f}\in\Ab_v$ generates the maximal ideal $\bar{\m}:=\m/\p_v$.

Assuming otherwise, there would exist an element $\bar{h}\in\bar{\m}$ such that $\bar{h}^2\mid \bar{f}$. Lift $\bar{h}$ to a polynomial $h\in K[x]$. Then $v_\xi(h)>v(h)=0$. By \cite[Theorem 5.1]{MacLane36} this implies $g\mid_v h$, which in turn implies
\[
    g^2\mid_v f.
\]
By the definition of $f$ and $u_v$ it follows that
\[
   g\mid_v \phi_1^{b_1}\cdot\ldots\cdot \phi_n^{b_n}, \qquad
    b_i := -e_vv(g)a_i\geq 0.
\]
But as a key polynomial, $g$ is $v$-irreducible (Definition \ref{def:key_polynomial}). We conclude that $g\mid_v \phi_i$, for some $i$. But this would mean that the prime ideal $\p_{\phi_i}\lhd A_v$ contains the maximal ideal $\m$. We have already argued in the proof of (i) that this is not the case, if we assume either Condition (b) or (c). This finishes the proof of the lemma.
\end{proof}

Let us now assume that $z$ lies on exactly two vertical components, $z\in E_{v_1}\cap E_{v_2}$, with $v_1<v_2$. Then by Proposition \ref{prop:component_tree}, $v_1,v_2$ are adjacent, i.e.~there is no other element of $V$ strictly in between $v_1$ and $v_2$. Together with Assumption \ref{ass:V} (c) this implies  $v_0:=\pred(v_2)\leq v_1$. We may therefore write
\[
    v_1 = [v_0,v_1(\phi)=\lambda_1], \quad
    v_2 = [v_0,v_2(\phi)=\lambda_2],
\]
where $\phi$ is a key polynomial for $v_0$ and $\lambda_1<\lambda_2\in\QQ$ (note that the case $v_0=v_1$ is not excluded). Let $N:=e_{v_0}$. Write $\lambda_i=b_i/c_i$, with $\gcd(b_i,c_i)=1$.

\begin{lem} \label{lem:regularity_cond2}
  Assume that
  \[
      \lambda_2 -\lambda_ 1 = \frac{N}{\lcm(N,c_1)\lcm(N,c_2)}.
  \]
  Then $z$ is a regular point of $X$ in which the components $E_{v_1}$ and $E_{v_2}$ intersect transversally. In other words: $D$ is a normal crossing divisor at $z$.
\end{lem}

\begin{proof}
This is \cite[Lemma 7.4]{ObusWewers}.
\end{proof}

\subsection{The algorithm} \label{subsec:algorithm}

Let $D_K\subset X_K=\PP^1_K$ be a reduced and effective divisor. We will now formulate an algorithm for constructing a regular model $X$ of $X_K$ such that $D:=D^\hor\cup X_s$ is a normal crossing divisor. See Example \ref{exa:running_example_rnc} for an explicit example worked out in detail.

Let $V_\infty$ denote the set of infinite pseudovaluations corresponding to the points of $D_K$. For simplicity, we assume that $V_\infty\subset V(K[x])^*$ consists of MacLane pseudovaluations. This is equivalent to saying that $D_K$ is the zero locus of a monic, integral and separable polynomial $f\in K[x]$.

\begin{algorithm} \label{alg:regular_model1}
We enlarge the set $V_\infty$ in three steps to oversets
\[
    V_\infty\subset V_1^*\subset V_2^*\subset V_3^*=V^*,
\]
as follows.

\begin{enumerate}[(1)]
\item
  $V_1^*$ is the union of $V_\infty$ with the set of all predecessors of elements of $V_\infty$,
  \[
     V_1^* := V_\infty \cup
       \bigcup_{v\in V_\infty} P(v).
  \]
  See Definition \ref{def:predecessors}.
\item
  The set $V_2^*$ is obtained by computing the {\em inf-closure} of $V_1^*$, i.e.\ the smallest subset of $V(K[x])^*$ containing $V_1^*$ and being closed under taking the infimum of two elements. In fact, it is not hard to show that
  \[
     V_2^* = V_1^* \cup \{\inf(v,w) \mid v,w\in V_1^*\}.
  \]
\item
  The set $V_3^*$ is obtained by applying Algorithm \ref{alg:regular_model2} below to every $v\in V_2:=V_2^*\cap V(K[x])$.
\end{enumerate}
\end{algorithm}

\begin{rem} \label{rem:algorithm1}
  The set $V_2^*$ produced by Step 2 of Algorithm \ref{alg:regular_model1} satisfies Assumption \ref{ass:V}. For Condition (a) this is true by construction. By Step 1, $V_1^*$ satisfies (c). Then in Step 2 we have added, a finite number of times, a new element which is less than some old element. By Lemma \ref{lem:predecessors}, this does not destroy Property (c). At the end of Step 2, $V_2^*$ satisfies (b) by construction.
\end{rem}

In Algorithm \ref{alg:regular_model2} below, we use shortest $N$-paths between rational numbers as defined in \cite{obus2022explicit}.

\begin{defn}
	If $a>a'\geq0$ are rational numbers, and $N$ is a positive integer, an {\em $N$-path} from $a$ to $a'$ is a sequence $a = b_0/c_0>b_1/c_1> \dots > b_n/c_n = a'$ of rational numbers in lowest terms such that
	\[
	\frac{b_i}{c_i} - \frac{b_{i+1}}{c_{i+1}} = \frac{N}{\lcm(N,c_i)\lcm(N,c_{i+1})}
	\]
	for $0\leq i\leq n-1$. If, in addition, no proper subsequence $b_0/c_0 > \dots > b_n/c_n$ containing $b_0/c_0$ and $b_n/c_n$ is an $N$-path, then the sequence is called a {\em shortest $N$-path} from $a$ to $a'$.
\end{defn}

We remark that for any pair of rational numbers $a>a'\geq0$ and any positive integer $N$, there exists a unique shortest $N$-path from $a$ to $a'$, \cite[Proposition A.14]{obus2022explicit}.

\begin{algorithm} \label{alg:regular_model2}
  Let $V^*$ be a finite set of MacLane pseudovaluations satisfying Assumption \ref{ass:V}. Let $v\in V:=V^*\cap V(K[x])$. We enlarge $V^*$ to an overset $V_v^*$, as follows.
  \begin{enumerate}[(1)]
  \item
    Let $v'\in V^*$ be such that $v<v'$ and that there is no element of $V^*$ strictly between $v$ and $v'$. Set $v_0:=\pred(v')$. By Assumption \ref{ass:V} (c), $v_0\leq v$, and we can write
    \[
        v = [v_0, v(\phi)=\lambda], \quad
        v' = [v_0,v'(\phi)=\lambda'].
    \]
    Set $N:=e_{v_0}$. Then we add to $V^*$ the MacLane valuations
    \[
       v_t := [v_0,v_t(\phi)=t],
    \]
    where $\lambda<t<\lambda'$ runs through the {\em shortest $N$-path from $\lambda'$ to $\lambda$}.

    We do this for all $v'$ as above.
  \item
    Let $v_0:=\pred(v)$ and write $v = [v_0,v(\phi)=\lambda]$.
    Set $N:=e_{v_0}$. If $\lambda\not\in N^{-1}\ZZ$ and $v'(\phi)=\lambda$ for all $v'$ as in (1), then we add to $V^*$ the MacLane valuations
    \[
       v_t := [v_0,v_t(\phi)=t],
    \]
    where $\lambda<t\leq \lambda'$ runs through the {\em shortest $N$-path from $\lambda'$ to $\lambda$}, and where
    \[
        \lambda' := \min \{ t\in \frac{1}{N}\ZZ \mid
           t>\lambda \}.
    \]
  \end{enumerate}
\end{algorithm}

\begin{rem} \label{rem:algorithm2}
  After applying Step (3) of Algorithm \ref{alg:regular_model1}, it is easy to see that the set $V^*=V_3^*$ still satisfies Assumption \ref{ass:V}.
\end{rem}

\begin{thm} \label{thm:regular_model_of_proejctive_line}
  Let $D_K\subset X_K$ be as at the beginning of this subsection. Let $V^*$ be the set of MacLane pseudovaluations produced by Algorithm \ref{alg:regular_model1} and let $(X,D^\hor)$ be the model of $(X_K,D_K)$ corresponding to $V^*$ via Proposition \ref{prop:models_valuations}. Then $X$ is regular, and $D:=D^\hor\cup X_s$ is a normal crossing divisor.
\end{thm}

\begin{proof}
By Remark \ref{rem:algorithm2}, the set $V^*$ satisfies Assumption \ref{ass:V}. Therefore, we only have to verify that the conditions of Lemma \ref{lem:regularity_cond1} and Lemma \ref{lem:regularity_cond2} are satisfied, for each closed point $z\in D^{\rm vert}$.

Consider a closed point $z$ on a vertical component $E_w$, $w\in V$. As we have shown in \S \ref{subsec:regularity}, we may assume that $\phi_w(z)\neq\infty_w$. This implies that $z$ either does not lie on any other vertical component (the situation of Lemma \ref{lem:regularity_cond1}), or $z\in E_{w'}$ where $w'>w$ is a MacLane valuation adjacent to $w$ (the situation of Lemma \ref{lem:regularity_cond2}). In both cases, the application of Algorithm \ref{alg:regular_model2} to the set $V^*$ guarantees that
\[
   w=v_t = [v_0,v(\phi)=t]
\]
where $t$ is a value from a `shortest path' from $\lambda'$ to $\lambda$ (we use the notation of Algorithm \ref{alg:regular_model2}).

If Condition (a), (b) or (c) of Lemma \ref{lem:regularity_cond1} hold then $z$ is a regular point and $D$ is a normal crossing divisor at $z$. We may therefore assume that all these three Conditions are false. This excludes the possibility that $t=\lambda'$, where $\lambda'$ is as in Part (2) of Algorithm \ref{alg:regular_model2}. Indeed, $t=\lambda'\in 1/N\ZZ$ implies that Condition (b) of  Lemma \ref{lem:regularity_cond1} is true, contrary to our assumption. We may therefore assume that $t<\lambda$.

Further, since Condition (c) is not satisfied, we may assume that the point $z\in E_w$ corresponds to the residue class
\[
    D_w(\m) = \{ w' \mid w'(\phi)>t \}.
\]
Let $t'>t$ be the next value after $t$ in the shortest path (this could be $\lambda'$). Then $w':=v_{t'}\in D_w(\m)\cap V$. This means that $z$ also lies on the vertical component $E_{w'}$. But by the definition of `shortest path', the pair $w=v_t$, $w'=v_{t'}$ satisfies the condition from Lemma \ref{lem:regularity_cond2}. Then this Lemma says that $z$ is a regular point and $D$ a normal crossing divisor at $z$. This concludes the proof of the theorem.
\end{proof}

\begin{rem} \label{rem:algorithm4}
  The algorithm described above has been implemented in the Sage/Python module\\ {\tt models\_of\_projective\_line} which is part of the Sage/Python package {\tt regular\_models}, \cite{regular_models}.
\end{rem}

%% file: global_sections.tex

\section{Computing the lattice of integral differential forms} \label{sec:global_sections}

Let $Y_K$ be a superelliptic curve, i.e.\ the smooth projective $K$-model of an affine plane curve given by an equation 
\[
  y^n = f(x).
\]
We assume that $n$ is invertible in the ring $\oo_K$ (Assumption \ref{ass:tame}). In this section, we present an algorithm for computing 
the integral differential forms (Definition \ref{def:integral_differential_forms}) of the curve $Y_K$.

Denote by $D_K \subset \PP_K^1$ the divisor defined by the polynomial $f$. Using Algorithm \ref{alg:regular_model1}, we construct a regular model $X$ of $X_K=\PP^1_K$ such that  $D:=D^{\hor} \cup X_s$ is a normal crossing divisor. By Theorem \ref{thm:tame}, the normalization of this model in the function field of $Y_K$ is a model $Y$ of $Y_K$ with only rational singularities; it can therefore be used for computing the lattice of integral differential forms. In this section we explain in detail how this is done. 
 
First, we briefly explain how to find, under some extra assumption, a $K$-basis of 
\[
   M_K = H^0(Y_K, \Omega_{Y_K/K})
\]
for a smooth projective curve $Y_K$. This is followed by a general discussion of reduced bases, which allows us to compute bases of integral differential forms for models that correspond to exactly one valuation on the function field $F_Y$. For the explicit computations of such reduced bases, it is necessary to compute the order of vanishing of a rational section along components of the special fiber $Y_s$. For our situation, where the model $Y$ is given as a cover of a model of the projective line, these computations are explained in the following subsection. Finally a conjunction of these results leads to Algorithm \ref{alg:superellipitccurves}.

\subsection{A $K$-basis}
\label{subsec:MK}

The first step is to find a $K$-basis of the vector space
\[
   M_K := H^0(Y_K,\Omega_{Y_K/K}).
\]
It is well known how to do this, but the precise algorithm depends of course on the way the curve $Y_K$ is defined. Let us assume, for simplicity, that $Y_K$ is defined generically by an equation $F(x,y)=0$, where $F\in K[x,y]$ is an absolutely irreducible polynomial, and the intersection of the plane model $F(x,y)=0$ of $Y$ with the torus $\mathbb{G}^2\subset\AA^2_K$ is smooth. Then there is a basis for $H^0(Y_K,\Omega_{Y_K/K})$ of the form
\[
x^{i-1}y^{j-1} \frac{dx}{F_y}, \quad (i,j) \in I,
\]
where $I\subset \ZZ^2$ are the points in the interior of the Newton polygon of $F$. See e.g.\ \cite[\S 1]{dokchitser_regular_models}. 

\begin{exa} \label{exa:basis_newton}
	Let $Y_K$ be superelliptic, with equation
	\[
	Y_K:\; y^n = f(x),
	\]
	with $n\geq 2$ and $f$ of degree $d\geq 3$. The condition that this equation defines a smooth curve in $\mathbb{G}^2_K$ means that $f$ has no multiple factors, except possibly an arbitrary power $x^m$ of $x$. A basis for $M_K$ is then 
	\[
	\left( x^{i-1} y^{j-n} dx \right)_{(i,j) \in I}, \quad I = \{(i,j): m(n-j) < ni < d(n-j)\} \subset \NN^2.
	\]
\end{exa}

In the following, we assume that a $K$-basis of $M_K$ is already known. 
Let $V(Y)$ denote the finite set of discrete valuations on the function field $F_Y$ of $Y$ corresponding to the vertical components of $Y$. 
Since $\omega_{Y/S}$ is a divisorial sheaf, 
we have
\[
   M = \{ \omega\in M_K \mid v(\omega) \geq 0 \;\forall\, v\in V(Y)\},
\]
see \eqref{eq:superelliptic2}.
Here $v(\omega)$ denotes the order of vanishing of $\omega$ along the component corresponding to $v$. See \eqref{eq:order_of_vanishing} for a precise definition. 

We choose a nonvanishing rational section $\eta$ of $\omega_{Y/S}$ (typically, $\eta=dx$) and obtain a $K$-linear embedding 
\[
M_K \inj F, \quad \omega \mapsto \omega/\eta.
\]
If we regard $M_K$ as a subvector space of $F_Y$, the definition of $M$ reads
\[
M = \{ f\in M \mid v(f) \geq -m_v \;\forall\, v\in V(Y) \},
\]
with $m_v:=v(\eta)$.

\subsection{Reduced bases}

It will be useful to consider a more abstract situation. 
Let $K$ be as above and $M_K$ a finite-dimensional $K$-vector space. Let $F/K$ be a field extension and $v:F\to\frac{1}{e_v}\ZZ\cup\{\infty\}$ a discrete valuation whose restriction to $K$ is equal to $v_K$. Here $e_v$ denotes the ramification index of the extension of valuations $v/v_K$. Let us also fix a $K$-linear embedding $M_K \inj F$ and an element $m\in \frac{1}{e_v}\ZZ$. The intermediate  goal is to determine the $\oo_K$-submodule
\[
    M_{v,m}:= \{ f\in M_K \mid v(f)\geq -m \}.
\]

\begin{defn}
  A system $(f_1,\ldots,f_i)$ of elements of $M_K$ is called {\em reduced} (with respect to $v$) if for all linear combinations $f=\sum_j a_jf_j$, $a_j\in k$, we have
  \[
      v(f) = \min_j v(a_jf_j).    
  \]
\end{defn}

\begin{lem} \label{lem:RR1}
  There exists a reduced basis $(f_1,\ldots,f_n)$ of $M_K$.
\end{lem}

\begin{proof}
Starting with an arbitrary $K$-basis $(g_1,\ldots,g_n)$ of $M_K$, a reduced basis can be constructed inductively. 

Suppose that $(f_1, \ldots, f_m)$ is a reduced basis for $\langle g_1, \ldots ,g_m\rangle$. For the construction of an element $f_{m+1}$ such that $(f_1,\ldots, f_{m+1})$ is a reduced basis for $\langle g_1, \ldots ,g_{m+1}\rangle$, we may assume that $f_{m+1} = g_{m+1} + \sum_{i=1}^{m}a_i f_i$ for some $a_i \in K$. It is easy to see that we get a reduced basis if and only if $f_{m+1}$ is the element with maximal valuation of this form.

We set $h = g_{m+1}$. If $h$ does not have maximal valuation, there exists an element $h' = h + \sum_{i=1}^{m}a_i f_i$ so that $h \sim_v - \sum_{i=1}^{m}a_i f_i$. (see Definition \ref{def:key_polynomial}.(i)). 
For determining $a_1, \ldots, a_m$ with this property (if they exist), we may as well do computations in the residue field of $v$. For that purpose, let $h_0 \in \langle f_1, \dots , f_m\rangle$ be such that $v(h) = v(h_0)$. If such an element $h_0$ does not exist, we cannot have $v(h') > v(h)$ and we are done, i.e. $f_{m+1} = h$.
 
Otherwise, we check whether  $h/h_0$ viewed as an element in the residue field of $v$, lies in the span of $(\pi^{e_i} f_i/h_0 \mid 1\leq i \leq m)$, where the integers $e_i$ are chosen minimally with the property $v(\pi^{e_i} f_i/h_0) \geq 0$. If this is not the case, $h$ already has maximal valuation. Otherwise, we may lift the relation that we obtain in the residue field and thereby construct an element $h'$ of the desired form with strictly larger valuation than $h$. Iterating this process, we will eventually find an element with maximal valuation. This is the next basis element $f_{m+1}$. 
\end{proof}

\begin{cor} \label{cor:RR1}
  For all $m\in\frac{1}{e_v}\ZZ$, the subset 
  \[
      M_{v,m} := \{ f\in M_K \mid v(f)\geq m \} \subset M_K
  \]
  is a free $\oo_K$-module of rank $n$. 
\end{cor}

\begin{proof}
Let $(f_1,\ldots,f_n)$ be a reduced $K$-basis of $M_K$. Then 
\[
    M_{v,m} = \gen{ \pi^{k_1}f_1,\ldots,\pi^{k_n}f_n },
\]
where
\[
    k_i := \left\lceil {m-v(f_i)} \right\rceil.
\]
\end{proof}

\begin{cor} \label{cor:RR2}
  Let $M_K\inj F$ be as before, with pairwise distinct discrete valuations $v_i:F\to\frac{1}{e_i}\ZZ\cup\{\infty\}$ and elements $m_i \in \frac{1}{e_i}\ZZ$ for $i \in \{1, \dots, r\}$. Then
  \[
     M := \{f\in M_K \mid v_i(f)\geq -m_i\;\forall\, i\}
  \]
  is a free $\oo_K$-module of rank $n$. 
\end{cor}

\begin{proof}
By definition,
\[
    M = \cap_i M_i, \quad M_i = \{f\in M_K \mid v_i(f)\geq -m_i \}.
\]
Since $M_i\subset M_K$ is a full lattice by Corollary \ref{cor:RR1}, $M\subset M_K$ is a full lattice, too.
\end{proof}

\begin{rem}
  The proofs of Lemma \ref{lem:RR1} and Corollary \ref{cor:RR1}-\ref{cor:RR2} can be easily turned into an algorithm to compute an $\oo_K$-basis of $M$. 
  In our implementation \cite{regular_models}, this is done in the modules {\tt lattices} and {\tt RR\_spaces}.
\end{rem}

\subsection{Order of vanishing of a rational section}
In this subsection, we restrict our considerations to superelliptic curves. We assume that $Y_K$ is given by an affine equation $y^n = f(x)$, where $n$ is invertible in $\oo_K$. Moreover, we let $X$ be a regular model of the projective line such that $D = D^{\hor} \cup X_s$ is a normal crossing divisor, where $D^{\hor}$ is the horizontal divisor defined by $f$. Then the normalization of $X$ in the function field of $Y_K$ is a model $Y$ for $Y_K$. The valuations in the set $V(Y)$ correspond to extensions of the valuations in $V(X)$ to $F_Y$. The computation of such extensions is explained in \cite[\S 4.6.2]{RuethThesis}.

Let $\omega$ be a nonzero rational section of $\oomega_{X/S}$. For any $v \in V(X)$, we can write $\omega=g\cdot\omega_0$, where $g\in F_X$ and $\omega_0$ is a generator of $\oomega_{X/S}$ at the generic point of $E_v$. Recall that $E_v$ denotes the component of the special fiber $X_s$ which corresponds to the valuation $v$. We set
\begin{align}
\label{eq:order_of_vanishing}
v(\omega) := v(g) \in \frac{1}{e_v}\ZZ.
\end{align}
It is easy to see that this value is independent of the choice of $\omega_0$.
We say that $v(\omega)$ is the {\em order of vanishing} of $\omega$ along the component $E_v$.\footnote{Note that in our setting the order of vanishing is not necessarily an integer, since $v$ is not normalized. In the literature, the integer $v(\omega) \cdot e_v$ is usually called the order of vanishing. } The definition of $w(\omega)$ for a rational section of $\oomega_{Y/S}$ and a valuation $w \in V(Y)$ is completely analogous. 

Before we show how to compute the order of vanishing of a differential form  in our setting, we first consider a more general situation. 

\begin{prop}
	\label{prop:algebraic_nt}
	Let $M$ be a rational function field over $K$ and $L/M$ a finite and separable extension. Further let $v_L$ be a geometric valuation on $L$ with valuation ring $\oo_L$. Denote by $v_M$ the restriction of $v_L$ to $M$ and by $\oo_M$ the corresponding valuation ring. Let $\phi: \Spec \oo_L \to \Spec \oo_M$  be the induced morphism.
	
	Assume that $\widehat{\oo_L} = \widehat{\oo_M}[\alpha]$ for some $\alpha \in \widehat{\oo_L}$ and let $p_{\alpha}$ denote the minimal polynomial of $\alpha$ over $\widehat{\oo_M}$. Then for a differential form $\omega_0 \in \oomega_{\Spec \oo_M/\Spec \oo_K}$, it holds that
	\[
	v_L(\phi^*\omega_0) = v_L(p_{\alpha}'(\alpha)) + v_M(\omega_0) .
	\]
	In particular, If $L/M$ is tame and totally ramified of degree $e$, then $v_L(\phi^*\omega_0) =  (e-1)/e_L+ v_M(\omega_0)$.	
\end{prop}

\begin{proof}{
	First, recall that \[
	\oomega_{\Spec \oo_L/\Spec\oo_K} = \oomega_{\Spec\oo_L/\Spec\oo_M} \otimes_{\oo_L} \phi^* \oomega_{\Spec\oo_M/\Spec\oo_K},
	\]
	by the adjunction formula, \cite[Theorem 6.4.9]{liu2002algebraic}. In the following, we describe $\oomega_{\Spec\oo_L/\Spec\oo_M}$. Since \[
	\oomega_{\Spec \oo_L/\Spec \oo_M}
	= (\oomega_{\oo_L/\oo_M})^{\tilde{}},\]
	it suffices to study the $\oo_L$-module $\oomega_{\oo_L/\oo_M}$. The latter is isomorphic to the the codifferent, $\mathfrak{C}(\oo_L/\oo_M)$ of the extension $\oo_L/\oo_M$ (\cite[Corollary A.2]{morrow2010explicit}). Recall that the different ideal  $\mathfrak{D}(\oo_L/\oo_M)\, \triangleleft\, \oo_L$ of the extension is defined as the complement of $\mathfrak{C}(\oo_L/\oo_M)$, i.e. \break ${\mathfrak{C}(\oo_L/\oo_M)\,\mathfrak{D}(\oo_L/\oo_M)=\oo_L}$.  
	
	For the computation of the different, we may pass to the completions of $\oo_L$ and $\oo_M$. More precisely,  it holds that \[\mathfrak{D}(\oo_L/\oo_M)  \;\widehat{\oo_L} = \mathfrak{D}(\widehat{\oo_L}/\widehat{\oo_M}),\] \cite[Satz 2.2]{neukirch1992algebraische}. By assumption $\widehat{\oo_L} = \widehat{\oo_M}[\alpha]$. In that situation the different ideal is generated by $p_{\alpha}'(\alpha)$, \cite[Satz 2.4]{neukirch1992algebraische}, hence $\mathfrak{C}(\oo_L/\oo_M)/\oo_L$ has length $v_L(p_{\alpha}'(\alpha)) \cdot e_L$.
	
	Let $\omega_0 \in \oomega_{\Spec \oo_M/\Spec \oo_K}$. Clearly, $\phi^*\omega_0 = 1 \otimes_{\oo_L} \phi^*\omega_0$ in the representation given by the adjunction formula. The above discussion implies
	\[
	v_L(\phi^*\omega_0) =  v_L(p_{\alpha}'(\alpha)) + v_M(\omega_0).
	\]
	
	If the extension $L/M$ is tame and totally ramified, then we can choose uniformizers $u_L$ and $u_M$ such that $\widehat{\oo_M}[u_L] = \widehat{\oo_L}$ and $p_{u_L} = T^e-u_M$. This means $p_{u_L}'(u_L) = e \cdot u_L^{e-1}$ and the statement follows from the formula in the general case.
}\end{proof}

We will now show how to compute $v(dx)$ for a valuation $v \in V(K[x])$. From that one can then deduce the value of $v(\omega)$ for any rational section $\omega \in H^0(X_K, \oomega_{X_K/K})$.

\begin{defn}
	\label{def:defining_system}
	Let $v \in V(X)$.
	We say that $(u_v,t_v)$ is a {\em defining system for $E_v$} if $u_v,t_v \in F_X = K(x)$ and 
	\begin{enumerate}[(i)]
		\item $u_v$ is a uniformizer for $v$.
		\item $\overline{t}_v$ is a separable transcendental generator of $k(v)$, i.e. $k(v)/k(\overline{t}_v)$ is finite and separable.
	\end{enumerate}
\end{defn}

\begin{prop} \label{prop:lci}
	Let $(u_v,t_v)$ be a defining system for $E_v$. Then there exists an irreducible polynomial $F_v \in \widehat{\oo_K[t_v]}[T]$ of degree $e_v$ satisfying $F_v(u_v) =0$. Moreover 
	\[ v(dx) = v(F_v'(u_v)) - v({dt_v}/{dx}). \]
\end{prop}

\begin{proof}{
	Consider the extension of function fields $K(x)/M$ with $M = K(t_v)$. This extension factors in the following way. 
	\[
	\begin{tikzcd}
	K(x) \arrow[swap]{d}{\phi_2} & v \arrow{d}{\text{res}}\\
	L = K(t_v,u_v) \arrow[swap]{d}{\phi_1} & v_L = [v_0,v_1(u_v) = 1/e_v, v_L(F_v) = \infty] \arrow{d}{\text{res}}\\
	M = K(t_v) & v_M
	\end{tikzcd}
	\]
	On the right, we included a description of the restrictions of $v$ to the respective fields. Note that the residue field of the valuation  $v_M$ is equal to $k(\bar{t_v})$ by definition of $t_v$. This implies that $v_M$ is the Gauß valuation on $K(t_v)$, i.e. the minimal element of $V(K[t_v])$.
	
	Since $L = K(t_v,u_v)$ has transcendence degree one over $K$, there exists an irreducible primitive polynomial $0 \neq F \in \oo_K[T_1,T_2]$ satisfying $F(t_v,u_v) = 0$. In particular $L =M[T]/(F_1)$, where $F_1 = F(t_v,T)$. Write \[
	0 = F_1(u_v) = \sum_{i=0}^n a_i u_v^i, \quad \text{where } a_i = \sum_{j=0}^{m_i}a_{i,j}t_v^j \in M.
	\]
	Since $a_i \in \oo_K[t_v] \subset \oo_M$ for all $i$, the polynomial is integral. Moreover the leading coefficient has valuation zero. Assume otherwise, i.e. $v_M(a_n)>0$. Having that $v_L(u_v^n)>v_L(u_v^i)$ for all $i \neq n$ implies that $v_M(a_i) >0$ for all $i$. Since $v_M$ is the Gauß-valuation on $K(t_v)$, this shows $v_M(a_{i,j})>0$ for all $i,j$.  
	
	Viewed as a polynomial over $\hat{M}$ (the completion of $M$ with respect to $v_M$), $F_1$ might be reducible. The irreducible factors are in correspondence with the extensions of $v_0$ to $L$ (see \cite[\S 4.6.2]{RuethThesis}). Note that we may scale  $F_1$ to make it monic in order to apply this correspondence. Since $v_M(a_n) = 0$, the polynomial remains integral.
	Let us denote by $F_v$ the factor corresponding to the valuation $v_L$. It follows from the description of $v_L$ that $F_v$ has degree $e_v$. 
	
	In order to compute $v(dx)$, first note that $v_M(dt_v) = 0$.  It follows from  Proposition \ref{prop:algebraic_nt} that  $v_L(\phi_1^*dt_v) = v_L(F_v'(u_v))$. Recall that $(t_v,u_v)$ is a defining system for the valuation $v$, in particular $u_v$ is a uniformizer for $v$. This implies that the extension of valuations $v \mid v_L$ is unramified that is $v_{\mid_L} = v_L$. Applying Proposition \ref{prop:algebraic_nt} again, we find that $v(\phi_2^* dt_v) = v_L(dt_v)$. Further note that $t_v \in L \subset K(x)$, hence we may denote $dt_v = \phi_2^*dt_v$.   Finally  $v(dt_v) = v(dt_v/dx) + v(dx)$ by definition. This concludes the proof.
}
\end{proof}
\begin{rem}
	\label{rem:find_defining_system}
	There exist methods in the Computer Algebra System Sage to compute a uniformizing element, $u_v$, and a generator for the residue field, $\overline{t_v}$. Lifting the latter to the ring $\oo_K$, one obtains $t_v$. Both $u_v$ and $t_v$ may be chosen to lie in $K[x]$. For more details, we refer the reader to \cite{RuethThesis}.
	
	In order to find the algebraic relation between $u_v$ and $t_v$, set 
	\[I = (T_1-u_v, T_2-t_v) \in K[T_1,T_2,x].
	\]
	Then $F$ is given by a generator of the elimination ideal $I \cap K[T_1,T_2]$. Using the methods outlined in \cite[\S 4.7]{RuethThesis} it is possible to compute arbitrarily good approximations of the factors of $F$ over the completion.
\end{rem}

\begin{lem}
	\label{lem:differentials_ramified}
	Let $X$ and $Y$ be as defined in the beginning of this subsection. 
	Let $\omega \in H^0(X_K,\Omega_{X_K/K})$. Then for any $w \in V(Y)$, there exists a valuation $v \in V(X)$ such that $\phi(E_w) = E_v$. Moreover 
	
	\[
	w(\phi^*\omega) = v(\omega) + \frac{e - 1}{e_w},
	\]
	where $e$ is the ramification index of the extension of valuations $[w:v]$. 
\end{lem}

\begin{proof}{The first part of the statement is clear from the construction of $Y$. 	Setting 
	\[
	L = F_Y, \; v_L = w, \quad M =F_X, \; v_M = v,
	\]
	the  second part follows from Proposition \ref{prop:algebraic_nt}.
}  
\end{proof}

\subsection{An algorithm for Superelliptic Curves}	
Here, we describe an algorithm for computing a basis for the lattice of integral differential forms for a superelliptic curve. This algorithm has been implemented in the Sage/Python module {\tt superelliptic\_curves} which is part of the package {\tt regular\_models}, \cite{regular_models}.

\begin{algorithm} \label{alg:superellipitccurves}
	Let $Y_K$ be a superelliptic curve with equation
	\[
	Y_K:\; y^n = f(x),
	\]
	with $n\geq 2$, $f$ of degree $d\geq 3$ and assume that $f$ has no multiple factors, except possibly an arbitrary power $x^m$ of $x$.\footnote{This assumption is made for simplicity only and has no bearing on our results} Write $D_K \subset \PP_K^1$ for the divisor defined by $f$.
	\begin{enumerate}[(1)]
		\item Compute $V = V(X) \subset V(K[x])$, where $X$ is a regular model of $\PP_K^1$ with the property that $D = D^\hor\cup X_s$ is a normal crossing divisor. 
		\item Choose an element $\eta \in M_K =  H^0(Y_K,\oomega_{Y_K/K})$ and compute a basis $\mathcal{B}_0$ of $M_K$ viewed as a subspace of $F_Y$ under the embedding $\omega \mapsto \omega/\eta$.
		\item Let $W = \{w \in V(F_Y) \mid w \large\mid_{K[x]} ~\in V(X)\}$. For each valuation $w \in W$ compute a basis $\mathcal{B}_w$ for the module $M_w := \{g \in M_K \mid w(g) \geq - w(\eta)\}$.
		\item Compute an $\oo_K$- basis $\mathcal{B}$ for $M_Y = \cap_{w \in W} M_w$.
	\end{enumerate}
\end{algorithm}

\begin{rem}
	We briefly explain how to perform the computations in each step of the above algorithm. As in our implementation, we choose $\eta = dx/y^{n-1}$ in Step (2). This choice is based on the fact that for every regular differential form $\omega$, there exists a polynomial $g \in K[x,y]$ such that $\omega = g dx/y^{n-1}$. As a result, the output of the algorithm is more legible. Of course the following steps can be easily adapted to a different choice of $\eta$. 
	\begin{enumerate}[(1)]
		\item Apply Algorithm \ref{alg:regular_model1} to $D_K$.
		\item With $\eta = dx/y^{n-1}$, a basis for $M_K$ is given by \[ 
		\mathcal{B}_0 = \left(
		x^{i-1} y^{j-1} \mid (i,j) \in \NN^2, m(n-j) < ni < d(n-j) 
		\right),\]
		cf. Example \ref{exa:basis_newton}.
		\item This step is an application of the results from the previous subsection.
			\begin{enumerate}
			\item First, compute $v(dx)$. For that purpose, determine a defining system $(u_v,t_v)$ and (an approximation of) the polynomial $F_v \in \widehat{\oo_K[t_v]}[T]$ as in Proposition \ref{prop:lci}. It is explained in Remark \ref{rem:find_defining_system} how to find such elements algorithmically. 
			Proposition \ref{prop:lci} then allows to compute  $v(dx)$. More precisely, we have
			\[
			v(dx) = v(F_v'(u_v)) - v(dt_v/dx).
			\]
			\item Next, compute all extensions $W_v \subset V(F_Y)$ of $v$ to the function field $F_Y$. It is explained in \cite[\S 4.6.2]{RuethThesis} how to compute such extensions of valuations. Then $W = \cup_{v \in V(X)} W_v$.  As a consequence of Lemma \ref{lem:differentials_ramified}, we get $w_v(dx) = v(dx) + (e -1)/e_w$ for each extension $w_v \in W_v$, where $e$ is the ramification index $[w_v:v]$. So 
			\[w_v(\eta) = v(dx) - (n-1)\cdot w_v(y) +(e-1)/e_w.\]
			\item In order to compute a basis $\mathcal{B}_v$ for $M_v := \{g \in M_K \mid w_v(g) \geq -w_v(\eta) \}$, start with the basis $\mathcal{B}_0$ of $M_K$ and modify it as described in the proof of Lemma \ref{lem:RR1} to get a basis which is reduced with respect to $v$. Then proceed as in the proof of Corollary \ref{cor:RR1}.
		\end{enumerate}
		\item We have to compute the intersection of finitely many $\oo_K$-lattices. For two lattices $M_1,M_2\subset M_K$, we have
		\[
		    M_1\cap M_2 \cong \ker\Big(
		       M_1\oplus M_2 \to M_K, \quad (m_1,m_2)\mapsto m_1-m_2,
		                      \Big)
		\]
		A basis of this module can easily be computed using the Smith normal form. The general case is done by iteration.
	\end{enumerate}
\end{rem}

\begin{thm}
	Let $Y_K$ be a superelliptic curve, with equation
	\[
	Y_K:\; y^n = f(x),
	\]
	with $n\geq 2$, $f$ of degree $d\geq 3$ and assume that $f$ has no multiple factors, except possibly an arbitrary power $x^m$ of $x$. Let $\mathcal{B}$ be the system produced by Algorithm \ref{alg:superellipitccurves} and $\eta \in H^0(Y_K,\oomega_{X_K/K})$ be the rational section chosen in Step (2) of the algorithm, then 
	\[
	\left(b \cdot \eta \mid b \in \mathcal{B} \right)
	\]
	is a basis for the lattice of integral differential forms for $Y_K$.
\end{thm}

\begin{proof}{
Let $X$ be the model from Step (1) and $Y$ the normalization of $X$ inside $F_Y$. It follows from Theorem \ref{thm:tame}  that $Y$ is a model of $Y_K$ with at most rational singularities. By Corollary \ref{cor:lattice_rationalsingularities}, $M := H^0(Y, \oomega_{Y/S})$ is the lattice of integral differential forms for $Y_K$.

As outlined in the end of \S \ref{subsec:MK}, 
\[
M
= \{g \in  M_K \mid w(g) \geq -w(\eta) \;\; \forall w \in V(Y)\}.
\]
The elements in  $V(Y)$ correspond to the extensions of elements in $V(X)$ to $F_Y$. More precisely $V(Y) = \{w \mid  w \textrm{ is an extension of some  } v \in V \textrm{ to } F_Y\}$. 
}
\end{proof}

%% file: examples_new.tex

\section{Examples} \label{sec:examples}
	
In this section we show how to compute a basis for the lattice of integral differential forms for superelliptic curves of the form  
\[
Y: y^n = (x^2-5)^3-5^5 \quad \text{over } \QQ_5.
\]
Using Algorithm \ref{alg:superellipitccurves}, we can find a basis for $M_Y$ for any $n$ not divisible by $5$. 
The first step in the algorithm is the same for all values of $n$. This computation is explained in Example \ref{exa:running_example_rnc}.
Following this example, we explain the remaining computations for the cases $n=2$ and $n=3$ in Examples \ref{exa:hyperell} and \ref{exa:n3} respectively. 

The examples presented in this section, along with more examples, can be found in the Jupyter Notebook {\tt Examples.ipynb} in our GitHub repository \cite{regular_models}.

\begin{exa}
	\label{exa:running_example_rnc}
	 We consider the polynomial  \[
	 f = (x^2-5)^3-5^5 \in \QQ_5[x].
	 \] Let $D_K \subset \PP_{K}^1$ be the divisor of zeroes of $f$ joined by the point $\infty$. This corresponds to the branch divisor of the cover $Y \to \PP_{K}^1$ defined by $(x,y) \mapsto (x)$.\footnote{If $n \mid 6$, then $\infty$ is not in the branch locus of the cover. Since we do not construct minimal models, there is no harm in adding it to the divisor $D_K$.}
	 
	 We will now compute the set $V^*$ as described in Algorithm \ref{alg:regular_model1}. Since $f$ is irreducible, the set $V_{\infty}$ consists of only two valuations,
	\[
	V_{\infty} = \{v_{\infty}, v_f\},
	\]
	{where } $v_{\infty} = [v_0, v(x) = -\infty], \; v_f = [v_0,v_1(x) = 1/2, v_2(x^2-5) = 5/3, v_3(f) = \infty]$.
	In the first step we add the predecessors of the valuations in $V_{\infty}$ and obtain
	\[
	V_1^* = V_{\infty} \cup \{v_0, [v_0,v_1(x) = 1/2],[v_1, v_2(x^2-5) = 5/3]\}.
	\]
	This set is already inf-closed, so there is nothing to do in the second step, i.e. $V_2^* = V_1^*$. 
	
	Applying Algorithm \ref{alg:regular_model2} to the set $V_2$, we find
	\[
	V_3^* = V_2^* \cup \left\{v_1' = [v_0,v_1(x) = 1]\} \cup \{v_{2,\lambda} = [v_1, v_2(x^2-5) = \lambda] \mid \lambda \in \{3/2, 7/4,2\}\right\}.
	\] 
	
	\input{models_exa_1}
	
	Let $(X,D)$ be the model corresponding to $V^* = V_3^*$, and denote by $Y$ the normalization of $X$ in the function field of $Y$. By Theorems \ref{thm:regular_model_of_proejctive_line} and \ref{thm:tame},  $Y$ is a model for $Y$ with at most tame cyclic quotient singularities. 	
\end{exa}

\begin{exa}
	\label{exa:hyperell}
	Let $Y$ be the genus-$2$ hyperelliptic curve defined by
	\[
	Y:y^2 = (x^2-5)^3-5^5 \quad \text{over } \QQ_5.
	\]
	In Example \ref{exa:running_example_rnc}, we have already computed a model $\X$ of the projective line as is needed for Step (1) in Algorithm \ref{alg:superellipitccurves}. 
	
	In Step (2), we choose $\eta = dx/y$ and since $Y$ is a genus-$2$ curve, a basis for $M_K$ is
	\[
	\mathcal{B}_0 = (1,\;x),
	\] 
	where $M_K$ is viewed as a subspace of $F_Y$. 
	
	Next, we have to compute the extensions of the valuations in $V = V^* \cap V(K[x])$ to $F_Y$. We denote this set by $W$. For example the unique extension of the valuation 
	\[
	v_2 = [v_0, v_1(x) = 1/2, v_2(x^2-5) = 5/3] \in V
	\] 
	to $F_Y$ is given by 
	\[
	w = [v_2, w(y^2-(x^2-5)^3-5^5) = \infty].
	\]
	In order to compute the module $M_w$ for this valuation, we first compute $v_2(dx)$. A defining system for $v_2$ is given by 
	\[
	(u_v,t_v) = \left(\frac{(x^2-5)\cdot x}{5^2}, \frac{(x^2-5)^3}{5^5}\right).
	\] 
	These two elements satisfy an algebraic relation $F_v(u_v,t_v) = 0$, where $F_v$ is a degree-$6$ polynomial. Computing this polynomial and applying Proposition \ref{prop:lci}, one obtains $v(dx) = 2$. More details are given in \cite[Example 7.10]{kunzweiler2021thesis}.
	
	Since the extension of valuations $w/v$ is unramified, we conclude that $w(dx) = 2$ as well (Lemma \ref{lem:differentials_ramified}). Moreover $w(y^2) = v_2((x^2-5)^3+5^5) =5$, hence $w(\eta) = w(dx)-w(y) = -1/2$. 
	
	It is easy to see that $\mathcal{B}_0$ is already reduced with respect to $w$. And Corollary \ref{cor:RR1} shows that 
	\[
	\mathcal{B}_w = (5,x)
	\]  
	is a basis for $M_w$. 
	
	\begin{table}[h!]
		\centering
		$\begin{array}{l|ccccccc}
		\hline
		v & v_0 &v_1&v_{1'}&v_{2,3/2}&v_2&v_{2,7/4}&v_{2,2} \\
		\hline \hline
		v(dx)    & 0 & 1    & 1 & 3/2  & 2    & 2 & 2 \\
		w_v(\eta)& 0 & -1/2 & 0 & -1/2 & -1/2 & -1/2 & -1/2 \\
		\hline
		\end{array}$
		\caption{Data for Step (3) with $n=2$}
		\label{table:running_exa_orders}
	\end{table}
	
	The data for the remaining valuations is summarized in Table \ref{table:running_exa_orders}. The notation for the valuations in $V$ is the same as in Example \ref{exa:running_example_rnc}. Note that some of the valuations $v \in V$ have more than one extension to $F_Y$. However, the value of $v(dx)$ (respectively $w_v(\eta)$) is the same for these extensions. Therefore we do not distinguish between different extensions in our table. In the implementation, it is necessary to perform the computation for all extensions, because the reduced bases might differ.

	In total, we obtain the basis
	\[
	\mathcal{B} = (5,x)
	\]
	for $M_Y$ (viewed as a subspace of $F_Y$).
\end{exa}

\begin{exa}
	\label{exa:n3}
	Now, we consider the genus-$4$ curve birationally defined by 
	\[
	Y:y^3 = (x^2-5)^3-5^5 \quad \text{over } \QQ_5.
	\]
	Clearly, the first step of the algorithm yields precisely the same output as in the previous example. In Step (2), we choose $\eta = dx/y^2$ and the basis
	\[
	\mathcal{B}_0 = (1,x,x^2,y)
	\]
	of $M_K$. 
	
	Next, the extensions of the valuations in $V = V(\X)$ to $F_Y$ are computed. The relevant data is given in Table 
	\ref{table:running_exa_n3}.\footnote{Although a cover of degree three is considered, we often obtain two valuations lying above a valuation $v \in V$ (see the second row of the table). Recall that we do not assume that the residue field is algebraically closed. Extending the  field $\QQ_p$ by adding a primitive third root of unity, the $2$'s in the second row of Table \ref{table:running_exa_n3} become $3$'s, whereas the remaining values do not change, since the extension is unramified.}
	
	\begin{table}[h!]
		\centering
		$\begin{array}{l|ccccccc}
		\hline
		v        & v_0 &v_1&v_{1'}&v_{2,3/2}&v_2&v_{2,7/4}&v_{2,2} \\
		\hline \hline
		v(dx)    & 0 &  1   & 1   &  3/2    & 2  & 2 & 2 \\
		\# \{w \in V(F): w{\big|}_{ K(x)} = v\}
		& 2 & 2    & 2   & 2      & 1   & 1 & 1\\
		e(w:v)   & 1 & 1    & 1   & 1      & 1   & 3 & 3 \\
		w_v(\eta)& 0 & -1   & -1  & -3/2   & -4/3& -7/6 & -1 \\
		\hline
		\end{array}$
		\caption{Data for Step (3) with $n=3$}
		\label{table:running_exa_n3}
	\end{table}
	In this case, the reduced bases become more interesting. For example consider the valuation $v_2 = [v_0, v_1(x) = 1/2, v_2(x^2-5) = 5/3]$ and its unique extension $w$ to $F_Y$. The basis
	\(
	(1,x,x^2-5,y)
	\)
	is reduced with respect to $w$.
	From the table, we can deduce that 
	\[
	\mathcal{B}_w = (5^2, 5x,x^2-5, y)
	\] 
	is a basis for $M_w$.
	
	Repeating this computation for all valuations $w \in W$ and then computing the intersection of the resulting lattices, we obtain 
	\[
	\mathcal{B} = (5^2, \, 5x,\, x^2-5,\,y).
	\] 	
\end{exa}

%% file: models_exa_1.tex
\begin{figure}[h!]
\centering
\subfloat[$V_\infty$ \label{Step0}]{
	\begin{tikzpicture}[grow =right,edge from parent/.style = {}, level/.style={sibling distance = .7cm/#1,level distance =.1cm}] 
	\node {}
	child{node  [gray node,label=right:{$v_{f}$}]{}}
	child{node [gray node, label=right:{$v_{\infty}$}]{}};
	\end{tikzpicture}
}	
\qquad
\subfloat[$V_1^* =V_2^*$ \label{Step1}]{
\begin{tikzpicture}[, grow =right,-,level/.style={sibling distance = 1cm/#1,
	level distance =1.5cm}] 
\node  [solid node,label=above:{$v_0$}]{}

	child{ node  [solid node,label=above:{$v_1$}]{}
		child{ node  [solid node,label=above:{$v_{2}$}]{}
			child{ node  [gray node,label=right:{$v_{f}$}]{}
			}
		}
	}
	child{ node [gray node, label=right:{$v_{\infty}$}]{}
	}
; 
\end{tikzpicture}
}
\qquad

\subfloat[$V_3^*$ \label{Step3}]{
	\begin{tikzpicture}[, grow =right,-,level/.style={sibling distance = .7cm, level distance =1.5cm}, level 1/.style={sibling distance = 1cm,	level distance =1.5cm}] 
	\node  [solid node,label=above:{$v_0$}]{}
	child{ node [solid node, label=above:{$v_1$}]{}
		child{ node [solid node, label=right:{$v_1'$}]{}	
		}
		child{ node  [solid node,label=above:{$v_{2,3/2}$}]{}
			child{ node  [solid node,label=above:{$v_{2}$}]{}
				child{ node  [gray node,label=right:{$v_{f}$}]{}
				}
				child{ node  [solid node,label=above:{$v_{2, 7/4}$}]{}
					child{ node  [solid node,label=above:{$v_{2,2}$}]{}
					}
				}
			}
		}
	}	
	child{node [gray node, label=right:{$v_{\infty}$}]{}
	}
	; 
	\end{tikzpicture}
}
	\caption{Different steps in Algorithm \ref{alg:regular_model1} with input $f = (x^2-5)^3 - 5^5$}
\end{figure}
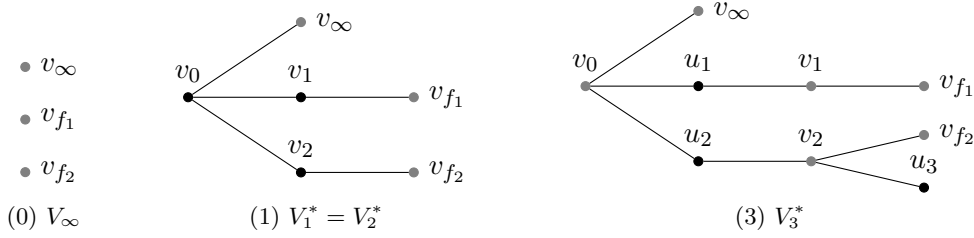